\documentclass[12pt, eqno]{article}
\usepackage{amssymb,latexsym}
\usepackage{amsmath,latexsym}
\usepackage{amscd,latexsym} 
\usepackage[all]{xy}

\newcommand{\N}{{\mathbb N}}

\newcommand{\R}{{\mathbb R}}

\newcommand{\Z}{{\mathbb Z}}

\newtheorem{theorem}{Theorem}[section]

\newtheorem{corollary}[theorem]{Corollary}

\newtheorem{definition}[theorem]{Definition}

\newtheorem{remark}[theorem]{Remark}

\newtheorem{lemma}[theorem]{Lemma}

\newtheorem{proposition}[theorem]{Proposition}

\begin{document}

\title{Some critical point
theorems and applications }

\author{\\Guangcun Lu
\thanks{Partially supported by the NNSF  10671017 and 10971014 of
China, and PCSIRT and Research Fund for the Doctoral Program Higher
Education of China (Grant No. 200800270003).}\\
{\normalsize School of Mathematical Sciences, Beijing Normal University},\\
{\normalsize Laboratory of Mathematics
 and Complex Systems,  Ministry of
  Education},\\
  {\normalsize Beijing 100875, The People's Republic
 of China}\\
{\normalsize (gclu@bnu.edu.cn)}}
\date{February 9,  2011} \maketitle \vspace{-0.1in}

\abstract{This paper is a continuation of \cite{Lu1}. In Part I,
applying the new splitting theorems developed therein we generalize
previous some results on computations of critical groups and some
critical point theorems to weaker versions. In Part II (in
progress), they are used to study multiple solutions for  nonlinear
higher order elliptic equations described in the introduction of
\cite{Lu1}.} \vspace{-0.1in}

\tableofcontents

\part{Some critical point theorems}

\section{Introduction }

In previous many critical point theorems involving computations of
critical groups the functionals are often assumed to be at least
$C^2$ smooth so that the usual splitting lemma can be used.
Doubtlessly, it is possible to obtain some new critical point
theorems or to generalize previous ones by combing our splitting
lemmas for continuously directional differentiable functionals with
and techniques and results in nonsmooth and continuous critical
point theories. Firstly,  we shall  generalize the some results in
\cite{BaLi, HiLiWa} to weaker versions.  Though they are not the
weakest versions, our theorems are more convenient in applications
because we do not need to compute such as subdifferentials and weak
slopes, which are not easy actually.  These are the main context in
Section 5. Next, in Section 6 we shall present the corresponding
version of the results on critical groups of sign-changing critical
points in \cite{BaChWa} and \cite{LiuWaWe} in our framework and
sketch how to prove them with our results in the previous sections.

\section{Compactness conditions and deformation
lemmas}\label{sec:2}
 \setcounter{equation}{0}

 Let $X$ be a normed vector space
with dual space $X^\ast$, $U\subset X$ nonempty and open, and
$f:U\to\R$ be a continuous functional. Recall in \cite{DegMa, IoSch,
Ka} that
 the {\bf weak slope}  of $f$ at $u\in U$ is the nonnegative extended real number
$|df|(u)$, which is   the supremum of the $\sigma$'s in $[0,
+\infty[$ such that there exist $\delta>0$ and ${\cal H}: B_X(u,
\delta)\times [0, \delta]\to X$ continuous with
 $$
 \|{\cal H}(v, t)-v\|\le t\quad\hbox{and}\quad f({\cal H}(v, t))\le f(v)-\sigma t.$$
 Clearly, $u\to |df|(u)$ is lower semicontinuous (\cite[Prop.2.6]{DegMa}).
 A point $u\in X$
 satisfying $|df|(u)=0$ is called a {\bf lower critical point} of
 $f$, and call $c=f(u)$ a {\bf lower critical value} of
 $f$.
By \cite[Def.1.1]{DeMaTo}, the {\bf strong slope} of a continuous
function $f:X\to\R$ at $u\in U$ is defined by
$$
|\nabla f(u)|=\left\{\begin{array}{ll}
 0     &\hbox{if}\; u\;\hbox{is a local minimum of}\; f,\\
 \overline\lim_{v\to u}\frac{f(v)-f(u)}{\|v-u\|} & \hbox{otherwise}.
 \end{array}\right.
 $$
Then $|df|(u)\le |\nabla f(u)|$ for any $u\in U$ (see below
Definition 2.8 in \cite{DegMa}), and
\begin{equation}\label{e:2.1}
|\nabla f(u)|\le \|f'(u)\|
\end{equation}
provided that $f$ has F-derivative $f'(u)$ at $u\in U$.

For a locally Lipschitz continuous function $f:U\to\R$, by
\cite[page 27]{Cl} it has the (Clarke) generalized gradient  at
every $u\in U$,
$$
\partial f(u)=\{g\in X^\ast\,|\, g(h)\le
f^\circ(u,h)\;\forall h\in X\},
$$
which is  the subdifferential at $\theta$ of the convex function
$X\to\R,\;h\mapsto f^\circ(u, h)$, where
$$
f^\circ(u, h)=\overline\lim_{ w\to\theta, s\downarrow 0}\frac{f(u+
w+ sh)-f(u+ w)}{s}
$$
is the generalized directional derivative of $f$ at $u$ in the
direction $h$.  If $X$ is a Banach
 space  it was proved in \cite{DegMa} that
\begin{equation}\label{e:2.2}
|df|(u)\ge |\partial f|(u):=\min\{\|x^\ast\|: x^\ast\in\partial
f(u)\}\quad \forall u\in U.
\end{equation}
This inequality may be strict by Example 1.1 in Krist\'{a}ly's
thesis \cite{Kr}.

 By \cite[Prop.2.2.1]{Cl} or
\cite[Prop.3.2.4(iii)]{Schi},  the function $f:U\to\R$ is strictly
H-differentiable at $u_0\in U$ if and only if  $f$ is locally
Lipschitz  continuous around $x_0$ and strictly G-differentiable at
$u_0\in U$. In this case we have $\partial f(u_0)=\{f'(u_0)\}$ by
\cite[Prop.2.2.4]{Cl}.  Moreover, by \cite[Prop.(6)]{Ch3} the
set-valued mapping $u\to\partial f(u)$ is  weak* upper
semi-continuous at $u_0$ in the sense that for any $\epsilon>0,
\,v\in X$ there exists a $\delta>0$ such that $|\langle w-f'(u_0),
v\rangle |<\epsilon$ for each $w\in\partial f(u)$ with $\|u-
u_0\|<\delta$. By \cite[Prop.(7)]{Ch3}  the function $u\to\|\partial
f\|(u)$ is lower semi-continuous at $u_0$, i.e.
$\underline{\lim}_{u\to u_0}\|f'(u)\|\ge\|f'(u_0)\|$. These are
summarized into  the following proposition.

\begin{proposition}\label{prop:2.1}
 Let $U$ be a nonempty open set of a normed vector
space $X$ with dual space $X^\ast$, and let $f:U\to\R$ be strictly
H-differentiable at every $u\in U$. Then \\
{\bf (i)}  $\partial f(u)=\{f'(u)\}$ at any $u\in U$, the map $U\to
X^\ast,\; x\mapsto f'(x)$ is weak* continuous and  the function
$U\ni u\to\|f'(u)\|$ is lower
semi-continuous.\\
{\bf (ii)}  $|\nabla f(u)|=|df|(u)=\|f'(u)\|\;\forall u\in U$
provided that $f$ is F-differentiable at $u\in U$.
\end{proposition}

 (ii) comes from (\ref{e:2.1}) and (\ref{e:2.2}).
Clearly, Proposition~\ref{prop:2.1} holds if $f$ is $C^1$.

Let $X$ be a Banach space and let $f:X\to\R$ be a G-differential
functional. Denote by $K(f)=\{x\in X\,|\, f'(x)=0\}$. For $c\in \R$
let $K(f)_c=\{x\in X\,|\, f(x)=c,\;f'(x)=0\}$ and $f^c=\{x\in X\,|\,
f(x)\le c\}$.  If $f$ is only continuous we write $LK(f)_c=\{x\in
X\,|\, f(x)=c,\;|df|(x)=0\}$.

\begin{definition}\label{def:2.2}
 {\rm Let $X$ be a Banach space and let  $f:X\to\R$ be a strictly
 H-differentiable functional. For $c\in\R$ the usual Palais-Smale
compactness condition at the level $c$, or $(PS)_c$ for short, means
that every $(PS)_c$ sequence $\{x_n\}\subset X$, i.e.  satisfying
$f(x_n)\to c$ and $f'(x_n)\to 0$, has a convergent subsequence;
moreover according to Cerami \cite{Cer} we say that $f$ satisfies
the condition $(C)_c$ if every $(C)_c$ sequence $\{x_n\}\subset X$,
i.e.  such that $f(x_n)\to c$ and $(1+ \|x_n\|)\|f'(x_n)\|\to 0$,
has  a convergent subsequence.}
\end{definition}

 Clearly, the second condition is weaker
 than the first one. Note that the semi-continuity of the map
 $u\mapsto\|f'(u)\|$ (by Proposition~\ref{prop:2.1}) implies
the limit of a $(PS)_c$ or $(C)_c$ sequence  $\{x_n\}$ is in
$K(f)_c$. In particular, $K(f)_c$ is compact if $f$ satisfies
$(PS)_c$ or $(C)_c$. The condition $(C)_c$ has the following
equivalent form
(\cite[Definition 1.1]{BBenFor}):\\
(i) every bounded sequence $(x_n)$ in $X$ with $f(x_n)\to c$ and
$f'(x_n)\to 0$ has a convergent subsequence; \\
(ii) there exist positive constants $\sigma, R, \alpha$ such that
$\|f'(x)\|\cdot\|x\|\ge\alpha$ for any $x$ with $c-\sigma\le f(x)\le
c+\sigma$ and $\|x\|\ge R$.

\begin{lemma}\label{lem:2.3}{\rm (First Deformation Lemma)}
For a strictly  H-differentiable functional $f:X\to\R$ on a (real)
Banach space $X$,   and $c\in\R$, suppose that $f$ satisfies the
condition  $(C)_c$. Then for every $\varepsilon_0>0$, every
neighborhood $U$ of $K(f)_c$ (if $K(f)_c=\emptyset$ we take
$U=\emptyset$), there exist an $0<\varepsilon<\varepsilon_0$ and a
map $\eta\in C([0, 1]\times X, X)$  satisfying
\begin{description}
\item[(i)] $\|\eta(t,u)-u\|\le {\rm e}(1+ \|u\|) t$, where ${\rm e}=\sum^\infty_{n=0}\frac{1}{n!}$;

\item[(ii)] $\eta(t, x)=x$ if $x\notin f^{-1}([c-\varepsilon_0, c+ \varepsilon_0])$;

\item[(iii)] $\eta\bigl(\{1\}\times (f^{c+\varepsilon}\setminus U)\bigr)\subset f^{c-\varepsilon}$;

\item[(iv)] $f(\eta(s, x))\le f(\eta(t, x))$ if $s\ge t$;

\item[(v)] $\eta(t, x)\ne x \Longrightarrow f(\eta(t, x))< f(x)$.
\end{description}
In particular, (ii)-(iv) show that $f$ satisfies the {\bf
deformation condition} $(D)_c$ in the sense of \cite[Def.
3.1]{BaLi}. When $f$ is even, $\eta$ may be chosen so that
$\eta(t,\cdot)$ is odd for all $t\in [0,1]$.
\end{lemma}

When $f$ is $C^1$ and satisfies the $(PS)_c$ (resp. $(C)_c$) this
lemma was proved by Palais \cite{Pa0} (see also \cite{Ra}), (resp.
Cerami \cite{Cer} and Bartolo-Benci-Fortunato \cite{BBenFor}). For a
$C^{1-0}$-functional on a  reflexive Banach space $X$, when $f$
satisfies the $(PS)_c$ (resp. $(C)_c$)  Chang \cite{Ch3} (resp.
Kourogenis and Papageorgiou\cite{KoPa}) proved this lemma.\

\noindent{\bf Proof of Lemma~\ref{lem:2.3}}.
 The ideas are following those of \cite[Theorem 4]{KoPa}.
  In the present case the
reflexivity of $X$ is not required because we do not need to use the
Eberlein separation theorem as in the proofs of \cite[Lemma
3.3]{Ch3} and \cite[Lemma 3]{KoPa}. Let us reprove Lemma 3 of
\cite{KoPa} under our assumptions as follows. By \cite[Lemma
2]{KoPa},  for each $\delta>0$ there exist $\gamma>0, \varepsilon>0$
such that for $K_c=K(f)_c$,
$$
(1+ \|x\|)\|f'(x)\|\ge \gamma \quad\forall x\in
f^{-1}([c-\varepsilon, c+ \varepsilon])\setminus N_\delta(K_c),
$$
where $N_\delta(K_c)=\{x\in X\,|\, d(x, K_c)<\delta\}$. For each
$x\in f^{-1}([c-\varepsilon, c+ \varepsilon])\setminus
N_\delta(K_c)$, we have $\|f'(x)\|\ge \gamma/(1+ \|x\|)$. Note that
$\|f'(x)\|=\sup\{\langle f'(x), h\rangle\,|\, h\in X,\,\|h\|=1\}$.
We have $h_x\in X$ such that $\|h_x\|=1$ and
$$
\langle f'(x), h_x\rangle>\frac{3\gamma}{4(1+
\|x\|)}>\frac{\gamma}{2(1+ \|x\|)}.
$$
By Proposition~\ref{prop:2.1}(i) we have $r_x>0$ such that
$$
\langle f'(y), h_x\rangle>\frac{\gamma}{2(1+ \|x\|)}\quad\forall
y\in B_X(x, r_x).
$$
 Now $\{B_X(x, r_x)\}$ is an open cover of $f^{-1}([c-\varepsilon,
c+ \varepsilon])\setminus N_\delta(K_c)$. Repeating the remaining
arguments in the proof of \cite[Lemma 2]{KoPa} we get a locally
Lipschitz vector field $V: f^{-1}([c-\varepsilon, c+
\varepsilon])\setminus N_\delta(K_c)\to X$ such that
$$\|V(x)\|\le (1+ \|x\|)\quad\hbox{and}\quad \langle f'(x),
V(x)\rangle\ge\frac{\gamma}{2}.
$$

Shrinking $\varepsilon>0, \gamma>0, \delta>0$ so that
$N_{3\delta}(K_c)\subset U$ and that (8) of \cite{KoPa} is
satisfied, and  almost repeating the proof of \cite[Theorem 4]{KoPa}
we may get the desired conclusions. The unique point which should be
noted is the proof of (iii).   By contradiction, suppose that
$f(\eta(1,x))>c-\varepsilon$ for some $x\in
f^{c+\varepsilon}\setminus U$. Then $c-\varepsilon<f(\eta(t,x))\le
c+\varepsilon$ for all $t\in [0, 1]$. As in the proof (c) of
\cite[Theorem 4]{KoPa} it must hold that $\eta([0,
1]\times\{x\})\cap (K(f)_c)_{2\delta}=\emptyset$. It follows that
there exist $0\le t_1<t_2\le 1$ such that $d(\eta(t_1,x), K(f)_c)=
2\delta$, $d(\eta(t_1,x), K(f)_c)= 3\delta$ and $2\delta<d(\eta(t,
x), K(f)_c)<3\delta$ for all $t_1\le t\le t_2$. Repeating the
remaining part of the proof (c) of \cite[Theorem 4]{KoPa} yields
(iii).  $\Box$\vspace{2mm}

Corresponding to \cite[Corollary 3.3]{BaLi} we have

\begin{corollary}\label{cor:2.4}
For a strictly  H-differentiable functional $f:X\to\R$ on a (real)
Banach space $X$, we have
\begin{description}
\item[(i)] If $f$ satisfies
the condition  $(C)_c$ for all $c\in [a, b]$ and $K(f)_c=\emptyset$
for $c\in [a, b]$ then there exist a deformation $\eta_t:X\to X$
such that $\eta_0=id_X$, $\eta_t(x)=x$ if $x\notin f^{-1}([a-1,
b+1])$, $f(\eta_t(x))$ is decreasing in $t$ and $\eta_1(f^b)\subset
f^a$.

\item[(ii)] If $f$ satisfies the condition $(C)_c$ for all $c\ge a$ and $K(f)_c=\emptyset$
for $c\ge a$ then there exist a deformation $\eta_t:X\to X$ such
that $\eta_0=id_X$, $\eta_t(x)=x$ if $f(x)\le a-1$, $f(\eta_t(x))$
is decreasing in $t$ and $\eta_1(X)\subset f^a$.
\end{description}
\end{corollary}

\noindent{\bf Proof}. We only outline the proof of (i). For each
$c\in [a, b]$ Lemma~\ref{lem:2.3} yields positive numbers
$\varepsilon^{(c)}_1<\varepsilon^{(c)}_2<1$ and a deformation
$\eta^{(c)}_t:X\to X$ such that $\eta^{(c)}_0=id_X$ and\\
$\bullet$ $\eta^{(c)}(t, x)=x$ if $x\notin
f^{-1}([c-\varepsilon^{(c)}_2, c+
\varepsilon^{(2)}_2])$;\\
$\bullet$ $\eta^{(c)}\bigl(\{1\}\times (f^{c+\varepsilon^{(c)}_1}\bigr)\subset f^{c-\varepsilon^{(c)}_1}$;\\
$\bullet$ $f(\eta^{(c)}(s, x))\le f(\eta^{(c)}(t, x))$ if $s\ge
t$.\\
Since $[a, b]$ is compact there exist finite numbers $a\le
c_1<\cdots<c_k\le b$ such that $\{(c_i-\varepsilon^{(c_i)}, c_i+
\varepsilon^{(c_i)})\}^k_{i=1}$ is an open cover of $[a, b]$. In
particular we have
$$
c_1-\varepsilon_1^{(c_1)}<a\le c_1<\cdots<c_k\le
b<c_k+\varepsilon_1^{(c_k)}.
$$
Then the composition
$\eta_t=\eta^{(c_1)}_t\circ\cdots\circ\eta^{(c_k)}_t$ satisfies the
desired requirements. As in the proof of \cite[Corollary
3.3(b)]{BaLi} we can derive (ii) from (i). $\Box$\vspace{2mm}

\begin{remark}\label{rm:2.5}
{\rm Even if $f$ is only continuous, if we replace  the $(C)_c$
condition by the following $(PS)_c$ condition
$$
``f(x_n)\to c\quad\hbox{and}\quad |df|(x_n)\to 0\Longrightarrow
\exists\;\hbox{a convergence subsequence of}\;\{x_n\}"
$$
 then Lemma~\ref{lem:2.3}  holds
provided that $K(f)_c$ is replaced by $LK(f)_c$. See \cite[Theorem
2.14]{CorDeg}. If $f$ is $F$-differentiable so that $|\nabla
f(x_n)|=|df|(x_n)=\|f'(x_n)\|\;\forall n$, then under the $(PS)_c$
condition Lemma~\ref{lem:2.3} is a corollary of \cite[Theorem
2.14]{CorDeg}. By the same reason, we may get the part I of the
following the second deformation lemma from  Theorem~2.3 of
\cite{Cor} or Theorem 4 and Remark 2 of \cite{Cor1}.}
\end{remark}

\begin{lemma}\label{lem:2.6}{\rm (Second Deformation Lemma)}
For a F-differentiable functional $f:X\to\R$ on a (real) Banach
space $X$, and $-\infty<a<b\le +\infty$ suppose that $f$ has only a
finite number of critical points at the level $a$ and has no
critical values in $(a,b)$.  Then\\
\noindent{\bf I}. If $f$   satisfies the condition $(PS)$ on
$f^{-1}([a, c])$ for all $c\in [a, b]\cap\R$, then  there exists a
deformation $\eta: [0, 1]\times f^{b\circ}\to f^{b\circ}:=\{f<b\}$
such that
\begin{description}
\item[(a)] $f(\eta(t,u))\le f(u)$;

\item[(b)] $u\in K(f)_a\Longrightarrow\eta(t, u)=u$;

\item[(c)] $\eta(\{1\}\times f^{b\circ})\subset f^{a\circ}\cup K(f)_a$;

\item[(d)] if $b\in\R$ and $K(f)_b=\emptyset$, then $\eta$ can be extended to $[0, 1]\times X$, still denoted
by $\eta$, such that $\eta(\{1\}\times f^{b})\subset f^{a\circ}\cup
K(f)_a$.
\end{description}
In particular, $f^{a\circ}\cup K(f)_a$ is a weak deformation retract
of $f^{b\circ}$.\\
\noindent{\bf II}.  If $f$ is $C^1$ and satisfies the condition
$(C)_c$ for all $c\in [a, b]\cap \R$, then $f^a$ is a strong
deformation retract of $f^b\setminus K(f)_b$, i.e. there exists a
map $\eta: [0, 1]\times (f^b\setminus K(f)_b)\to (f^b\setminus
K(f)_b)$, called a strong deformation retraction of $f^b\setminus
K(f)_b$ onto $f^b$, satisfying
\begin{description}
\item[(i)] $\eta(0,u)=u$ for all $u\in f^b\setminus K(f)_b$;

\item[(ii)] $\eta(t, u)=x$ for all $(t, u)\in [0, 1]\times f^a$;

\item[(iii)] $\eta(\{1\}\times (f^{b}\setminus f^a))=f^a$.
\end{description}
\end{lemma}

For the part II, under the condition $(PS)_c$ the proof is due to
Rothe \cite{Ro}, Chang \cite{Ch4} and Wang \cite{Wa0}; and under the
condition $(C)_c$ the proof can be found in Bartsch-Li \cite{BaLi},
Perera-Schechter \cite{PerSche} and
  Perera-Agarwal-O'Regan \cite{PerAO}.

Applying these two deformation lemmas and our splitting lemma,
Theorem 2.1 in \cite{Lu1}, the standard arguments as in \cite{Ch,
MaWi, Mi, Pa0} may yield the following two theorems.

\begin{theorem}\label{th:2.7}
Let $H$ be a Hilbert space and let $f:H\to \R$ be a
$F$-differentiable and strictly $H$-differentiable functional.
Suppose:
\begin{description}
\item[(i)] for some small $\varepsilon>0$ there exists a unique
critical value $c$ in  $[c-\varepsilon, c+\varepsilon]$,
\item[(ii)] $K_c$ is finite and $f$ satisfies the conditions of Theorem
2.1 in \cite{Lu1}  near each of $K_c$,

\item[(iii)] {\rm Either} $f$ satisfies the (PS) condition on
$f^{-1}([c-\varepsilon, c+\varepsilon])$ {\rm or} $f$ is $C^1$ and
satisfies the condition $(C)_d$ for every $d\in [c-\varepsilon,
c+\varepsilon]$.
\end{description}
Then for any abelian group $G$ one has
$$
H_\ast(f_{c+\varepsilon}, f_{c-\varepsilon};G)\cong\bigoplus_{z\in
K_c}C_\ast(f, z),
$$
which is finitely dimensional vector spaces over $G$ if $G$ is a
field.
\end{theorem}

Let $B^m$ be the closed unit disk in $\R^m$. By a {\it topological
embedding} $h:B^m\to H$ we mean that it is continuous bijection onto
$h(B^m)\subset H$ and that $h$ is a  homeomorphism between $B^m$ and
$h(B^m)$ with respect to the induced topology on $h(B^m)$ from $H$.

\begin{theorem}\label{th:2.8}
{\rm (Handle Body Theorem)}.\quad Under the assumptions of
Theorem~\ref{th:2.7}, if each of $K_c=\{z_j\}^l_1$ is also
nondegenerate, then for some $0<\epsilon\le\varepsilon$ there exist
topological embeddings  $h_i:B^{m_i}\to H$, $i=1,\cdots, l$, such
that
$$
f_{c-\epsilon}\cap h_j(B^{m_j})=f^{-1}(c-\epsilon)\cap
h_j(B^{m_j})=h_j(\partial B^{m_j})
$$
for $j=1,\cdots,l$, and
$f_{c-\epsilon}\cup\bigcup^l_{j=1}h_j(B^{m_j})$ is a deformation
retract of $f_{c+\epsilon}$, where $m_j$ is the Morse index of
$z_j$.
\end{theorem}

Similarly, we can also give the versions on Hilbert manifolds.

The following is a slight variant of \cite[Prop. 2.1]{HiLiWa}.

\begin{proposition}\label{prop:2.9}
Let $H$ be a Hilbert space and let $B(\infty):H\to H$ be a bounded
self-adjoint linear operator satisfying (${\bf C1_\infty}$), i.e.
$0$ is at most an isolated point of the spectrum
$\sigma(B(\infty))$, which implies $\pm(B(\infty)u, u)_H\ge
2\alpha_\infty\|u\|^2\;\forall u\in H^\pm_\infty$. Assume:
\begin{description}
\item[(i)] $g:H\to\R$ is  strictly
H-differentiable (i.e. locally Lipschitz continuous and strictly
G-differentiable) and hence $\partial g(x)=\{g'(x)\}$ by
Proposition~\ref{prop:2.1};
\item[(ii)] $\|g'(x)\|$ is bounded, $g'$ is compact and
$\nu_\infty=\dim H^0_\infty<\infty$;
\item[(iii)] For any $M>0$, $g'(u^0+ u^\pm)\to 0$ uniformly in $u^\pm\in\bar B_H(\theta, M)\cap H^\pm_\infty$
as $\|u^0\|\to\infty$.
\end{description}
Then ${\cal L}(u)=\frac{1}{2}(B(\infty)u,u)_H+  g(u)$ satisfies (PS)
condition on $H\setminus C_{R, M}$, where
$$
C_{R,M}=\{u=u^0+ u^\pm\,|\, \|u^0\|>R,\;\|u^\pm\|<M\}.
$$
Consequently, for any $(PS)_c$ sequence $\{u_n\}$ of ${\cal L}$,
{\rm either} $\{u_n\}$ has a bounded subsequence (and hence a
converging subsequence) {\rm or} $c\in C_\infty({\cal L})$ and there
exists a subsequence $\{u_{n_k}\}$ such that $\|u^0_{n_k}\|\to
\infty$, $\|u^\pm_{n_k}\|\to 0$ and $g(u_{n_k})\to c$. Here
$C_\infty({\cal L})$ is a closed subset of $\R$ given by
\begin{eqnarray*}
&&C_\infty({\cal L}):=\{c\in\R\;|\; \exists\, u^0_n\in
H^0_\infty,\;u^\pm_n\in H^\pm_\infty\quad\hbox{with}\quad
\|u^0_n\|\to\infty,\;\nonumber\\
&&\hspace{45mm}\|u^\pm_\infty\|\to 0\;\hbox{such that}\quad g(u^0_n+
u^\pm_n)\to c\}.
\end{eqnarray*}
Consequently,  ${\cal L}$ satisfies the $(PS)_c$ condition for
$c\notin C_\infty$.
\end{proposition}

\noindent{\bf Proof}. Let $\{u_n\}\subset H\setminus C_{R, M}$ be
such that  ${\cal L}(u_n)\to c$ and $B(\infty)u_n+ g'(u_n)\to 0$ as
$n\to\infty$. Since $\|g'(x)\|$ is bounded, and $\|u^\pm_n\|\le
\|B(\infty)|_{H^\pm_\infty}\|\cdot\|B(\infty)u^\pm_n\|$ we infer
that $\{\|u^\pm_n\|\}$ is bounded. Let $\|u^\pm_n\|\le M_1\;\forall
n$. Suppose that a subsequence $\|u^0_{n_k}\|\to\infty$. Then
$g'(u^0_{n_k}+ u^\pm_{n_k})\to 0$ and so
\begin{eqnarray*}
\|u^\pm_{n_k}\|&\le&
\|B(\infty)|_{H^\pm_\infty}\|\cdot\|B(\infty)u^\pm_{n_k}\|\\
&=& \|B(\infty)|_{H^\pm_\infty}\|\cdot\|{\cal
L}'(u_{n_k})-g'(u^0_{n_k}+ u^\pm_{n_k})\|\to 0.
\end{eqnarray*}
Hence $u_{n_k}=u^0_{n_k}+ u^\pm_{n_k}\in C_{R,M}$ for $k$ large
enough. This contradiction shows that $\{\|u^0_n\|\}$ is bounded.
Since $g'$ is compact and $\nu_\infty=\dim H^0_\infty<\infty$ we
have a subsequence $\{u_{n_k}\}$ such that $u^0_{n_k}\to u^0$ and
$g'(u_{n_k})\to v$. The latter implies that
\begin{eqnarray*}
\|u^\pm_{n_k}-u^\pm_{n_l}\|&\le&
\|B(\infty)|_{H^\pm_\infty}\|\cdot\|B(\infty)u^\pm_{n_k}-B(\infty)u^\pm_{n_l}\|\\
&=& \|B(\infty)|_{H^\pm_\infty}\|\cdot\|{\cal L}'(u_{n_k})-{\cal
L}'(u_{n_l})-[g'(u_{n_k}) -g'(u_{n_l})]\|\to 0
\end{eqnarray*}
as $k, l\to\infty$. Hence $\{u_{n_k}\}$ converges to some $v$.
$\Box$\vspace{2mm}

\section{Computations of critical groups}\label{sec:3}
 \setcounter{equation}{0}

\subsection{Critical groups at infinity and
computations}\label{sec:3.1}

 In this subsection ${\bf K}$ always denotes a commutative ring without
 special statements.  For a strictly H-differentiable
functional $f:X\to\R$ on a Banach space $X$, suppose that the set of
critical values of $f$ is strictly bounded from below by $a\in\R$,
and  for all $c\le a$ that $f$ satisfies  the condition $(C)_c$. By
Corollary~\ref{cor:2.4}(i),  for every nonnegative integer $m$,
\begin{eqnarray}
&&C_m(f,\infty; {\bf K}):=H_m(X, f^{a}; {\bf K}),\label{e:3.1}\\
&&C^m(f,\infty; {\bf K}):=H^m(X, f^{a}; {\bf K})\label{e:3.2}
\end{eqnarray}
are independent of the choices of such $a$, and are called the $m$th
critical group of $f$ at infinity and $m$th cohomological critical
group of $f$ at infinity, respectively (cf. Definition 3.4 of
\cite{BaLi}). Here $H_\ast(-;{\bf K})$ and $H^\ast(-;{\bf K})$
denote the singular homology and cohomology with coefficients in
${\bf K}$. It is well-known that
\begin{eqnarray*}
C_m(f,\infty;{\bf K})=H_m(X, f^{a}; {\bf K})&\cong &{\rm Hom}(H_m(X,
f^{a}; {\bf K}))\\
&\cong& H^m(X, f^{a}; {\bf K})=C^m(f,\infty;{\bf K})
\end{eqnarray*}
if ${\bf K}$ is a field. Let $\bar H^\ast(-;{\bf K})$ denote
Alexander-Spanier cohomology with coefficients in ${\bf K}$, which
has often some stronger excision  and continuity properties. Now the
Banach space $X$ is a ANR (absolute neighborhood retract). By
Section K on the page 30 of \cite{Hu1},  every open subset of an ANR
an ANR, and Hanner theorem claims that  a metrizable space is an ANR
if it has a countable open covers consisting of ANR. Hence
$f^a=\cup^\infty_{n=1}\{f<a-\frac{1}{n}\}$ is an ARN. From Section 9
of \cite[Chapter 6]{Spa} it follows that $H^m(X, f^{a}; {\bf
K})\cong \bar H^m(X, f^{a}; {\bf K})$ for any field ${\bf K}$. In
particular we have
\begin{equation}\label{e:3.3}
C_m(f,\infty;{\bf K})\cong C^m(f,\infty;{\bf K})\cong \bar H^m(X,
f^{a}; {\bf K})
\end{equation}
 for any field ${\bf K}$ and nonnegative integer $m$. These
and Proposition 3.15 of \cite{PerAO} lead to

\begin{proposition}\label{prop:3.1}
 For a strictly H-differentiable functional $f:X\to\R$ on a Banach
space $X$, suppose that the set of critical values of $f$ is
strictly bounded from below by $a\in\R$, and that $f$ satisfies  the
condition $(C)_c$ for all $c\in\R$. Then for any field ${\bf K}$ and
nonnegative integer $m$ it holds
\begin{description}
\item[(i)] $C_m(f,\infty;{\bf K})\cong C^m(f,\infty;{\bf
K})\cong\delta_{m0}{\bf K}$ if $f$ is bounded from below.

\item[(ii)] $C_m(f,\infty;{\bf K})\cong C^m(f,\infty;{\bf
K})\cong \widetilde H^{m-1}(f^a;{\bf K})$ if $f$ is unbounded from
below. Here $\widetilde H^0(f^a;{\bf K})= H^0(f^a;{\bf K})/{\bf K}$
and $\widetilde H^q(f^a;{\bf K})= H^q(f^a;{\bf K})$ for $q\ge 1$.
\end{description}
\end{proposition}

By Proposition~B.1 of \cite{Lu1},  the continuously directional
differentiability is stronger than the strict $H$-differentiability.
 We have the following generalization of Theorem~3.9 in \cite{BaLi}.

\begin{theorem}\label{th:3.2}
Suppose for $V_\infty=H$:
\begin{description}
\item[(i)]  the assumptions of Theorem~4.1 of \cite{Lu1}, (${\bf S}$), $({\bf
F1_\infty})$-$({\bf F3_\infty})$ and $({\bf C1_\infty})$-$({\bf
C2_\infty})$, $({\bf D_\infty})$ and $({\bf E'_\infty})$, are
satisfied;

\item[(ii)]   ${\cal L}(u)=\frac{1}{2}(B(\infty)u,u)_H+
o(\|u\|^2)$ as $\|u\|\to\infty$;

\item[(iii)] $\nabla{\cal L}(u)= B(\infty)u +
o(\|u\|)$ as $\|u\|\to\infty$, where $\nabla{\cal L}$ is the
gradient of ${\cal L}$ defined by $d{\cal L}(u)(v)=(\nabla{\cal
L}(u), v)_H$ for all $u, v\in H$; ({\rm Note}: we do not assume
${\cal L}\in C^1(H, \R)$.) \footnote{For a possible method removing
this condition, see below the end of this document.}

\item[(iv)] the critical values of ${\cal L}$ are bounded below;

\item[(v)] ${\cal L}$ satisfies  the condition $(C)_c$
(or $(D)_c$) for $c\ll 0$.
\end{description}
 Then $C_k({\cal L},\infty;{\bf K})=0$ for $k\in [\mu_\infty, \mu_\infty+
 \nu_\infty]$ even if $\mu_\infty=\infty$ or $\nu_\infty=\infty$. Moreover, if $\mu_\infty<\infty$
 and $\nu_\infty=0$ then $C_{\mu_\infty}({\cal L}, \infty;{\bf K})\cong{\bf K}$ (even if $H$ is not complete).
\end{theorem}

\noindent{\bf Proof.} {\it Step 1}. Carefully checking the proof of
Lemma 4.2 in \cite{BaLi} one easily sees that the conditions (i) and
(ii) imply: \textsf{for sufficiently large $R>0$ and $a\ll 0$ the
pair
$$
\bigl(B_{H^0_\infty}(\theta, R+1)\oplus H^\pm_\infty, {\cal L}^a\cap
(B_{H^0_\infty}(\theta, R+1)\oplus H^\pm_\infty)\bigr)
$$
is homotopy to the pair
$$
\bigl(B_{H^0_\infty}(\theta, R+1)\oplus \bar B_{H^-_\infty}(\theta,
1), B_{H^0_\infty}(\theta, R+1)\oplus \partial\bar
B_{H^-_\infty}(\theta, 1)\bigr).
$$
The homotopy equivalence leaves the $H^0_\infty$-component fixed. }

{\it Step 2}. Under the assumption (i), by Theorem~4.1 of \cite{Lu1}
we can get Lemma~4.3 of \cite{BaLi}:  \textsf{There exist a
sufficiently large $R>0$, $a\ll 0$ and a continuous map
$\gamma:B_{H^0_\infty}(\infty, R)\to [0, 1]$ with $\gamma(C)>0$ for
$C:=B_{H^0_\infty}(\theta, R+1)\cap B_{H^0_\infty}(\infty, R)$ such
that the pair
$$
(B_{H^0_\infty}(\infty, R)\times H^\pm_\infty, {\cal
L}^a\cap(B_{H^0_\infty}(\infty, R)\times H^\pm_\infty))
$$
is homotopy equivalent to the pair $(B_{H^0_\infty}(\infty, R)\times
H^-_\infty, \Gamma)$, where
\begin{eqnarray*}
&&\Gamma=\{(z, u)\in B_{H^0_\infty}(\infty, R)\times
H^-_\infty:\;\|u\|\ge\gamma(z)\}quad\hbox{and}\\
&& \gamma(z)=\left\{
\begin{array}{ll}
0 &\hbox{if}\; {\cal L}(z+ h^\infty(z))\le a,\\
1 &\hbox{if}\; {\cal L}(z+ h^\infty(z))\ge a+ 1,\\
{\cal L}(z+ h^\infty(z))-a &\hbox{elsewhere}.
\end{array}\right.
\end{eqnarray*}
Moreover, the  homotopy equivalence leaves the
$H^0_\infty$-component fixed.}

{\it Step 3}. By the assumptions (iv) and (v), $C_\ast({\cal
L},\infty;{\bf K})=H_\ast(H, {\cal L}^a;{\bf K})$ for $a\ll 0$ is
well-defined. Using Step 1 and Step 2 we may repeat the proof on the
pages 428-429 of \cite{BaLi} to obtain at the desired conclusion.
$\Box$\vspace{2mm}

Using Corollary~\ref{cor:2.4} we derive (i) and (ii) the following
proposition, which are corresponding with Propositions~3.5 and 3.6
in \cite{BaLi}.

\begin{proposition}\label{prop:3.3}
 For a strictly
H-differentiable functional $f:X\to\R$ on a Banach space $X$, we
have
\begin{description}
\item[(i)] If $a<\inf f(K(f))\le\sup f(K(f))<b$ and
$f$ satisfies the condition $(C)_c$ (or $(D)_c$) for any $c\notin
(a, b)$, then $C_\ast(f,\infty;{\bf K})\cong H_\ast(f^b, f^a;{\bf
K})$.

\item[(ii)] If $f$ satisfies  the condition $(C)_c$ (or $(D)_c$)  for any
$c\in \R$, then $C_\ast(f,\infty;{\bf K})\cong 0$ in the case
$K(f)=\emptyset$, and $C_\ast(f,\infty;{\bf K})\cong C_\ast(f,
x_0;{\bf K})$ in the case $K(f)=\{x_0\}$.

\item[(iii)] If $f$ is F-differentiable, satisfies the condition
$(PS)_c$ for every $c\in\R$ and has finite critical points, then for
every field ${\bf K}$ it holds that
$$
\dim C_m(f, \infty;{\bf K})\le\sum_{u\in K(f)}\dim C_m(f, u;{\bf
K})\quad\forall m\in\N\cup\{0\}.
$$
\end{description}
\end{proposition}

When $f$ is $C^1$ and satisfies the condition $(C)_c$ for every
$c\in\R$, (iii) was proved in Proposition 3.16 of \cite{PerAO}.

\noindent{\bf Proof of Proposition~\ref{prop:3.3}}. We only prove
(iii). It is almost standard (see \cite{Cor} and  \cite[Remark
 2]{Cor1}). For the reader's
convenience we prove it.  Since $f$ has only finite critical points
we may take $-\infty<a<\inf f(K(f))$. Let $c_1<\cdots<c_k$ be all
critical values. Take numbers $a<a_1<\cdots <a_k$ such that
$$
a_1<c_1<a_2<c_2<\cdots<a_k<c_k.
$$
Take $b=+\infty$ then $f^{b\circ}=X$, and $f^{a_1\circ}\subset
f^{a_2\circ}\subset\cdots\subset f^{a_k\circ}\subset f^{b\circ}$. By
\cite[Lemma 3.12]{PerAO} we get for every nonnegative integer $m$
that
$$
\dim\bar H^m(f^{b\circ}, f^{a_1\circ};{\bf K})\le
\sum^k_{i=1}\dim\bar H^m(f^{a_{i+1}\circ}, f^{a_i\circ};{\bf K}),
$$
where $a_{k+1}=b$. From the proof of Theorem~4 in \cite{Cor1} we may
see that $f^{a_1\circ}$ is a strong deformation retract of $f^b$.
This implies $\bar H^m(f^{b\circ}, f^{a_1\circ};{\bf K})\cong \bar
H^m(X, f^{a};{\bf K})$, and so
$$
\dim\bar H^m(f^{b\circ}, f^{a_1\circ};{\bf K})=\dim
C_m(f,\infty;{\bf K})\quad\forall m.
$$
Now the part I of Lemma~\ref{lem:2.6} yields
$$
\bar H^m(f^{a_{i+1}\circ}, f^{a_i\circ};{\bf K})\cong \bar
H^m(f^{a_{i+1}\circ}, f^{c_i\circ};{\bf K})\cong \bar
H^m(f^{c_i\circ}\cup K(f)_{c_i}, f^{c_i\circ};{\bf K})
$$
 for any $i=1,\cdots, k$ and nonnegative integer $m$. As showed in \cite[Remark
 2]{Cor1} the subsets $f^{c_i\circ}\cup K(f)_{c_i}$ and
 $f^{c_i\circ}$ are ARN. Hence
$$
 \bar H^m(f^{c_i\circ}\cup K(f)_{c_i}, f^{c_i\circ};{\bf K})\cong
H^m(f^{c_i\circ}\cup K(f)_{c_i}, f^{c_i\circ};{\bf K})\quad\forall
m, i.
$$
Let $K(f)_{c_i}=\{x_{i1},\cdots, x_{il_i}\}$, and $U_{i1},\cdots,
U_{il_i}$ be mutually disjoint open neighborhoods of $x_{i1},\cdots,
x_{il_i}$, respectively. Then the excision property of singular
cohomology groups lead to
$$
 H^m(f^{c_i\circ}\cup K(f)_{c_i}, f^{c_i\circ};{\bf K})\cong
\bigoplus^{l_i}_{j=1}H^m\bigl((f^{c_i\circ}\cup \{x_{ij}\})\cap
U_{ij}, f^{c_i\circ}\cap U_{ij};{\bf K})
$$
 for all $m, i$. Moreover, it is easy to prove
\begin{eqnarray*}
H^m\bigl((f^{c_i\circ}\cup \{x_{ij}\})\cap U_{ij}, f^{c_i\circ}\cap
U_{ij};{\bf K})&\cong& H^m\bigl(f^{c_i}\cap U_{ij},
(f^{c_i}\setminus\{x_{ij}\})\cap U_{ij};{\bf K})\\
&\cong& H_m\bigl(f^{c_i}\cap U_{ij},
(f^{c_i}\setminus\{x_{ij}\})\cap U_{ij};{\bf K})\\
&\cong& C_m(f, x_{ij};{\bf K}).
\end{eqnarray*}
The desirable conclusion follows from these immediately.
$\Box$\vspace{2mm}

We can also give  generalizations of some computation results on
critical groups at infinity  such as  Proposition 3.10 in
\cite{BaLi} and some parts of \cite{Li2}. We leave them intersecting
readers.

\subsection{Computations of critical groups at
degenerate critical points}\label{sec:3.2}

\begin{definition}\label{def:3.4}
 {\rm Let $X$ be a Banach space and let  $f:X\to\R$ be a strictly
 H-differentiable functional. For subsets $S\subset X$ and $A\subset X^\ast$ we call $f$ to satisfy
{\bf $(PS)$ condition with respect to $A$ on $S$} if every sequence
$(x_n)$ in $S$ with $(f(x_n))$ bounded and $f'(x_n)\to y\in A$ has a
convergent subsequence.}
\end{definition}

Clearly, the usual $(PS)$ condition is the $(PS)$ condition with
respect to $\{0\}$ on $X$. We have the following generalization of
Proposition 2.5 in \cite{BaLi}.

\begin{proposition}\label{prop:3.5}
\begin{description}
\item[(i)] Let $X\subset H$ be as in ({\bf S}) of \cite[\S 2.1]{Lu1}, and let ${\cal L}$
be a continuously directional differentiable functional defined in
an open neighborhood $V$ of $x_0\in X$ in $H$; moreover we assume
that the conditions ({\bf F1})-({\bf F3}), ({\bf C1})-({\bf C2}) and
({\bf D}) in \cite[\S 2.1]{Lu1} hold with $\theta$ replaced by
$x_0$.

\item[(ii)]  $x_0$ is an isolated critical point of ${\cal L}$, and  ${\cal L}$ is
$F$-differentiable near $x_0$.

\item[(iii)]  $\nabla{\cal L}(u)= B(x_0)(u-x_0) +
o(\|u-x_0\|)$ as $\|u-x_0\|\to 0$, where $\nabla{\cal L}$ is the
gradient of ${\cal L}$ defined by $d{\cal L}(u)(v)=(\nabla{\cal
L}(u), v)_H$ for all $u, v\in H$. ({\rm Note}: we do not assume
${\cal L}\in C^1(H, \R)$.)

\item[(iv)] ${\cal L}$  satisfies
 the (PS) condition with respect to $H^0$ on a closed ball $\bar B_H(x_0,
 \delta)$.
\end{description}
Let $\mu_0$ and $\nu_0$ be the Morse index and nullity of ${\cal L}$
at $x_0$. Then  we have:
\begin{description}
\item[(a)]  $C_k({\cal L}, x_0;{\bf K})=\delta_{k\mu_0}{\bf K}$ provided that ${\cal L}$ also
satisfies:\\
 {\bf (${\bf AC}^+$)} There exist $\varepsilon>0$ and $\theta\in (0,
 \pi/2)$ such that $(\nabla{\cal L}(u+ x_0), u^0)_H\ge 0$ for any $u=u^0+
 u^\pm\in H^0+ H^\pm$ with $\|u\|\le\varepsilon$ and
 $\|u^\pm\|\le\|u\|\cdot\sin\theta$.

\item[(b)]  $C_k({\cal L}, x_0;{\bf K})=\delta_{k(\mu_0+\nu_0)}{\bf K}$ provided that ${\cal L}$ also
satisfies:\\
 {\bf (${\bf AC}^-$)} There exist $\varepsilon>0$ and $\theta\in (0,
 \pi/2)$ such that $(\nabla{\cal L}(u+ x_0), u^0)_H\le 0$ for any $u=u^0+
 u^\pm\in H^0+ H^\pm$ with $\|u\|\le \varepsilon$ and
 $\|u^\pm\|\le\|u\|\cdot\sin\theta$.
\end{description}
 \end{proposition}

By Corollary~2.6 of \cite{Lu1} we have $C_q({\cal L}, x_0;{\bf
K})=0$ if $q\notin [\mu_0, \mu_0+ \nu_0]$. So
Proposition~\ref{prop:3.5} may be viewed a refinement of this
result. For the proof of it
 we also need the following stability property of critical groups for continuous functionals by
Cingolani and Degiovanni \cite{CiDe}, which is a very general
generalization of the previous results due to Chang \cite[page 53,
Th.5.6]{Ch}, Chang and Ghoussoub \cite{ChGh} and in Mawhin and
Willem \cite[Th.8.8]{MaWi}, and Corvellec and Hantoute \cite{CorHa}.

\begin{theorem}\label{th:3.6}{\rm (\cite[Th.3.6]{CiDe})}
Let $\{f_t:\,t\in [0, 1]\}$ be a family of continuous functions from
a metric space $X$ to $\R$, let $U$ be an open subset of $X$ and
$[0, 1]\ni t\mapsto u_t\in U$ a continuous map. Assume:
\begin{description}
\item[(I)] if $t_k\to t$ in $[0, 1]$, then $f_{t_k}\to f_t$
uniformly on $\overline{U}$;

\item[(II)] $\overline{U}$ is complete,  and for every sequence
$t_k\to t$ in $[0, 1]$ and $(v_k)$ in $\overline{U}$ with
$|df_{t_k}|(t_k)\to 0$ and $(f_{t_k}(v_k))$ bounded, there exists a
subsequence $(v_{k_j})$ convergent to some $v$ with $|df_t|(v)=0$;

\item[(III)] $|df_t|(v)>0$ for every $t\in [0, 1]$ and $v\in \overline{U}\setminus\{u_t\}$

\end{description}
Then $C_q(f_0, u_0;{\bf K})\cong C_q(f_1, u_1;{\bf K})$ for every
$q\ge 0$.
\end{theorem}

\noindent{\bf Proof of Proposition~\ref{prop:3.5}}. Following the
proof ideas of Proposition 2.5 in \cite{BaLi}, we assume
$x_0=\theta$. For the case (a) (resp. (b)) we set ${\cal
L}_t(u)={\cal L}(u)+ \frac{1}{2}t\|u^0\|^2$ (resp. ${\cal
L}_t(u)={\cal L}(u)- \frac{1}{2}t\|u^0\|^2$) for $t\in [0, 1]$.
Using the assumptions (ii) and (iv), as in the proof of Proposition
2.5 in \cite{BaLi} we have a small $\epsilon\in (0, 2/\delta)$ such
that $\theta$ is the only critical point of each ${\cal L}_t$ in
$B_H(\theta, 2\epsilon)$. Clearly, ${\cal L}_1$ also satisfies the
assumption of Theorem~2.1 of \cite{Lu1}, and $\theta$ is a
nondegenerate critical point of ${\cal L}_1$ with Morse index
$\mu_0$ (resp. $\mu_0+ \nu_0$) in the case (a) (resp. (b)). It
follows from Corollary~2.6 of \cite{Lu1} that
\begin{equation}\label{e:3.4}
C_k({\cal L}, \theta)=\delta_{k\mu_0}{\bf K}\;\hbox{in case
(a)}\quad(\hbox{resp.}\; C_k({\cal L},
\theta)=\delta_{k(\mu_0+\nu_0)}{\bf K}\;\hbox{in case (b)}).
\end{equation}
The remaining is to  prove $C_q({\cal L}, \theta)\cong C_q({\cal
L}_1, \theta)$ for every $q\ge 0$ in both cases.

We only prove the case (a). Since ${\cal L}$ is continuously
directional differentiable, and $F$-differentiable at $\theta$, so
is each ${\cal L}_t$. By Proposition~B.2(ii) of \cite{Lu1} and
Proposition~\ref{prop:2.1} every ${\cal L}_t$ is strictly
$H$-differentiable (and thus locally Lipschitz continuous), and
\begin{equation}\label{e:3.5}
\partial {\cal L}_t(u)=\{{\cal L}'_t(u)\}\quad\hbox{and}\quad  |d{\cal L}_t|(u)=\|{\cal
L}'_t(u)\|\;\forall u\in \bar B_H(\theta,\epsilon)
\end{equation}
and  $|d{\cal L}_t|(u)>0\;\forall u\in \bar
B_H(\theta,\epsilon)\setminus\{\theta\}$ because $\theta$ is only
critical point of ${\cal L}_t$ in $B_H(\theta, 2\epsilon)$.

The second equality in (\ref{e:3.5}) implies that $\theta$ is also a
lower critical point of each ${\cal L}_t$. Because of (\ref{e:3.4}),
we hope to use Theorem~\ref{th:3.6} proving that $C_q({\cal L},
\theta)\cong C_q({\cal L}_1, \theta)\;\forall q\ge 0$. It suffice to
check the conditions of Theorem~\ref{th:3.6}. Clearly, ${\cal
L}_{t_k}\to{\cal L}_t$ uniformly on $\bar B_H(\theta, \epsilon)$ as
$t_k\to t$ in $[0, 1]$. Now we assume: $t_k\to t$ in $[0, 1]$,
$(u_k)\subset \bar B_H(\theta, \epsilon)$ is such that $({\cal
L}_{t_k}(u_k))$ is bounded and $|d{\cal L}_{t_k}|(u_k)\to 0$. These
imply that ${\cal L}'_{t_k}(u_k)={\cal L}'(u_k)+ t_ku^0_k\to\theta$.
Since $\dim H^0<\infty$ we may assume $u^0_k\to u^0$ (passing a
subsequence if necessary). Then $({\cal L}(u_k))$ is bounded and
${\cal L}'(u_k)\to -tu^0\in H^0$. By the assumption (iv) $(u_k)$ has
a convergent subsequence $u_{k_i}\to u_0\in \bar B_H(\theta,
\epsilon)$. Hence $u^0=P^0u_0$. Moreover, since ${\cal L}$ is
continuously directional differentiable we get
$$
({\cal L}'(u_{k_i}), v)_H\to ({\cal L}'(u_0), v)_H\;\forall v\in
H.
$$
Hence $(-tu^0, v)_H=({\cal L}'(u_0), v)_H\;\forall v\in H$,
i.e. ${\cal L}'(u_0)+ tP^0u^0={\cal L}'(u_0)+ tu^0=\theta$. It
follows from (\ref{e:3.5}) that $|d{\cal L}_t|(u_0)=\|{\cal
L}_t(u_0)\|=0$. Namely $\{{\cal L}_t\,|\, t\in [0, 1]\}$ satisfies
the conditions of Theorem~\ref{th:3.6} on $\bar
B_H(\theta,\epsilon)$.
  $\Box$\vspace{2mm}

\section{Morse inequalities and some critical point theorems}\label{sec:4}
 \setcounter{equation}{0}

\subsection{Morse inequalities}\label{sec:4.1}

 In this subsection,
unless otherwise specified, let the functional ${\cal L}:H\to\R$ be
as in Proposition~\ref{prop:2.9}. Then ${\cal L}$ satisfies the
$(PS)_c$ condition for each $c\notin C_\infty({\cal L})$. We also
assume that ${\cal L}$ is $C^1$ so that it has   a pseudo-gradient
vector field, $V:\widetilde H\to H$. Note that $\widetilde H\to
H,\;u\mapsto V(u)/\|V(u)\|$ is also a locally Lipschitz continuous
map. Consider the flow
\begin{equation}\label{e:4.1}
\dot\eta(t,
u)=-\frac{V(\eta(t,u))}{\|V(\eta(t,u))\|}\quad\hbox{and}\quad
\eta(t,0)=u.
\end{equation}
Our goal is to present the Morse inequality established in
\cite{BaLi, HiLiWa} under the above weaker setting. For $F\subset H$
let $\widetilde F=\cup_{t\in\R}\eta(t, F)$.

 For any isolated value
$c$ in $C_\infty({\cal L})$, let ${\cal
L}^{c+\varepsilon}_{c-\varepsilon}:={\cal L}^{-1}([c-\varepsilon,
c+\varepsilon])$ and $K^{c+\varepsilon}_{c-\varepsilon}({\cal
L}):=K({\cal L})\cap {\cal L}^{c+\varepsilon}_{c-\varepsilon}$.
Define
\begin{eqnarray*}
&&U_{R, M}=\{u=u^0+ u^\pm\,|\, \|u^0\|\le R\}\cup\{u=u^0+ u^\pm\,|\,
\|u^0\|>R,\;\|u^\pm\|\ge M\},\\
&&C_{R, M}=\{u=u^0+ u^\pm\,|\, \|u^0\|> R,\,\|u^\pm\|< M\}=H\setminus U_{R, M},\\
&&U^{c+\varepsilon}_{R, M}=U_{R, M}\cap{\cal
L}^{c+\varepsilon},\quad
C^{c+\varepsilon}_{R, M}=C_{R, M}\cap{\cal L}^{c+\varepsilon},\\
&&A^{c+\varepsilon}_{R, M}=U^{c+\varepsilon}_{2R, M/2}\cap
C^{c+\varepsilon}_{R, M}.
\end{eqnarray*}

The following lemma corresponds to Lemma 2.1, Proposition 2.2 and
Corollary 2.4 in \cite{HiLiWa} (see also Lemma 2.7 and Theorem 2.9
in \cite{Li}).

\begin{lemma}\label{lem:4.1}
Assume that $K({\cal L})^{c+\varepsilon_0}_{c-\varepsilon_0}=K({\cal
L})_c$ is compact for some $\varepsilon_0>0$. Then for $R$ large and
$R>M>0$ with $K({\cal L})_c\subset B_H(\theta, R/2)\cup C_{3R, M/8}$
there exists $\varepsilon_1>0$ such that for any $\varepsilon\in
(0,\varepsilon_1)$ it holds that
\begin{description}
\item[(i)] $\bigl({\cal L}^{c+\varepsilon}_{c-\varepsilon}\cap \widetilde{U^{c+\varepsilon}_{R,
M}}\bigr)\cap \bigl({\cal L}^{c+\varepsilon}_{c-\varepsilon}\cap
\widetilde{C^{c+\varepsilon}_{2R, M/4}}\bigr)=\emptyset$,

\item[(ii)] ${\cal L}^{c+\varepsilon}\cap \widetilde{A^{c+\varepsilon}_{R,
M}}\cong {\cal L}^{c-\varepsilon}\cap
\widetilde{A^{c+\varepsilon}_{R, M}}$.
\end{description}
Furthermore, if $K({\cal L})_c$ is compact then for any $M>0$ there
exist a large $R>0$, and $\varepsilon_1>0$ such that for any
$\varepsilon\in (0,\varepsilon_1)$,
\begin{eqnarray*}
H_q({\cal L}^{c+\varepsilon}, {\cal L}^{c-\varepsilon}; {\bf
K})&\cong& H_q({\cal L}^{c+\varepsilon}\cap
\widetilde{U^{c+\varepsilon}_{2R, M/2}}, {\cal
L}^{c-\varepsilon}\cap
\widetilde{U^{c+\varepsilon}_{2R, M/2}};{\bf K})\\
&\oplus &H_q({\cal L}^{c+\varepsilon}\cap
\widetilde{C^{c+\varepsilon}_{R, M}}, {\cal L}^{c-\varepsilon}\cap
\widetilde{C^{c+\varepsilon}_{R, M}};{\bf K})\quad\forall
q=0,1,\cdots.
\end{eqnarray*}
Hereafter ${\bf K}$ always denotes a commutative ring without
 special statements.
\end{lemma}

\noindent{\bf Proof.} Since $U_{3R, M/8}\setminus B_H(\theta,
R/2)=H\setminus \bigl(B_H(\theta, R/2)\cup C_{3R, M/8}\bigr)$ is
disjoint with $K({\cal
L})^{c+\varepsilon_0}_{c-\varepsilon_0}=K({\cal L})_c$, and ${\cal
L}$ satisfies the $(PS)$ condition in $U_{3R, M/8}$ by
Proposition~\ref{prop:2.9}, there exists an $\varepsilon'>0$ such
that
$$
\|{\cal L}'(u)\|\ge\|{\cal L}'(u)\|\ge \varepsilon'\quad\forall u\in
{\cal L}^{c+\varepsilon_0}_{c-\varepsilon_0}\cap\bigl(U_{3R,
M/8}\setminus B_H(\theta, R/2)\bigr).
$$

For $0<\varepsilon<\min\{\varepsilon_0, \varepsilon'M/16\}$, suppose
that $\eta(s, x)\in {\cal L}^{c+\varepsilon}_{c-\varepsilon}\cap
C^{c+\varepsilon}_{R+ M/4, 3M/4}$ for some $x\in {\cal
L}^{c+\varepsilon}_{c-\varepsilon}\cap U^{c+\varepsilon}_{R, M}$ and
$s>0$. Then there exists $t_1<t_2$ such that
\begin{eqnarray*}
&&\eta(t_1, x)\in {\cal L}^{c+\varepsilon}_{c-\varepsilon}\cap
\partial U_{R, M},\quad
\eta(t_2, x)\in {\cal L}^{c+\varepsilon}_{c-\varepsilon}\cap
\partial C_{R+ M/4, 3M/4}\subset
{\cal L}^{c+\varepsilon}_{c-\varepsilon}\cap
\partial U_{2R, M/2},\\
&&\eta(t,x)\in \bar C_{R, M}\cap U_{R+ M/4, 3M/4}\quad\forall t\in
[t_1, t_2].
\end{eqnarray*}
Then $M/4\le\|\eta(t_1, x)-\eta(t_2,x)\|\le |t_2-t_1|$, and
\begin{eqnarray*}
{\cal L}(\eta(t_2, x))-{\cal
L}(\eta(t_1,x))&=&\int^{t_2}_{t_1}\frac{d}{dt}{\cal
L}(\eta(t,x))dt\\
&=&-\int^{t_2}_{t_1}\frac{\langle{\cal L}'(\eta(t,x)),
V(\eta(t,x))\rangle}{\|V(\eta(t,x))\|}dt\\
&\le &-\frac{\varepsilon'}{2}|t_2-t_1|.
\end{eqnarray*}
Hence
$$
\varepsilon'\frac{M}{4}\le\varepsilon'|t_2-t_1|\le 2\bigl({\cal
L}(\eta(t_1, x))-{\cal L}(\eta(t_2,x))\bigr)\le 4\varepsilon.
$$
This contradicts to the choice of $\varepsilon$. So we get
$$
\widetilde{\bigl({\cal L}^{c+\varepsilon}_{c-\varepsilon}\cap U_{R,
M}\bigr)}\cap \bigl({\cal L}^{c+\varepsilon}_{c-\varepsilon}\cap
C_{R+ M/4, 3M/4}\bigr)=\emptyset.
$$
Since $\bigl({\cal L}^{c+\varepsilon}_{c-\varepsilon}\cap
\widetilde{U^{c+\varepsilon}_{R, M}}\bigr)\subset
\widetilde{\bigl({\cal L}^{c+\varepsilon}_{c-\varepsilon}\cap U_{R,
M}\bigr)}$,
\begin{equation}\label{e:4.2}
\bigl({\cal L}^{c+\varepsilon}_{c-\varepsilon}\cap
\widetilde{U^{c+\varepsilon}_{R, M}}\bigr)\cap\bigl({\cal
L}^{c+\varepsilon}_{c-\varepsilon}\cap C_{R+ M/4,
3M/4}\bigr)=\emptyset.
\end{equation}
Similarly, for $0<\varepsilon<\min\{\varepsilon_0,
\varepsilon'M/16\}$ we have
$$
\bigl({\cal L}^{c+\varepsilon}_{c-\varepsilon}\cap
\widetilde{C^{c+\varepsilon}_{2R, M/4}}\bigr)\cap\bigl({\cal
L}^{c+\varepsilon}_{c-\varepsilon}\cap U_{2R- M/4,
M/2}\bigr)=\emptyset.
$$
This and (\ref{e:4.2}) together give (i).

As in the proof of \cite{HiLiWa, Li} we can get the remaining
conclusions. $\Box$\vspace{2mm}

Following \cite{BaLi, HiLiWa},  as in \cite{Ch} and \cite{MaWi})
Lemma~\ref{lem:4.1} may lead to

\begin{theorem}\label{th:4.2}
Under the assumptions of Proposition~\ref{prop:2.9}, let ${\cal L}$
be $C^1$, and let $K({\cal L})$ and $C_\infty({\cal L})$ be finite.
Denote by $\beta_k(f, x):=\dim C_k(f, x;{\bf K})$ for $x\in K({\cal
L})$, and by
\begin{eqnarray*}
&&\beta_k({\cal L}, c)=\dim H_k({\cal
L}^{c+\varepsilon}\cap\widetilde C_{R, M}, {\cal
L}^{c-\varepsilon}\cap\widetilde C_{R, M};{\bf K})\quad\forall c\in
C_\infty({\cal L}),\\
&&P({\cal L},\infty):=\sum^\infty_{k=0}\dim H_k(H, {\cal L}^a;{\bf
K})t^k
\end{eqnarray*}
for any $a<\min\bigl\{x\,|\, x\in {\cal L}(K({\cal L}))\cup
C_\infty({\cal L})\cup\{0\}\bigr\}$. Then there exists a polynomial
$Q(t)$ with nonnegative integer coefficients such that
$$
P({\cal L}, \infty)+ (1+ t)Q(t)=\sum_{x\in K(f)}P({\cal L}, x)+
\sum_{v\in C_\infty(f)}P({\cal L}, c),
$$
where $P({\cal L},x):=\sum^\infty_{k=0}\beta_k({\cal L}, x)t^k$ and
$P({\cal L},c):=\sum^\infty_{k=0}\beta_k({\cal L}, c)t^k$
\end{theorem}

\subsection{Some critical point theorems}\label{sec:4.2}

Many critical point theorems, which were obtained by computations of
critical groups, can be generalized with our methods. For example,
the following is a generalization of Theorem 5.1 on the page 121 of
\cite{Ch}.

\begin{theorem}\label{th:4.3}
{\rm (${\bf I}$)} Let the Banach space  $(X, \|\cdot\|)$ and the
Hilbert space $(H, (\cdot,\cdot)_H)$ satisfy the condition ({\bf S})
in \cite[\S 2.1]{Lu1}. Let $\widetilde H$ (resp. $\widetilde X$) be
a $C^1$ Hilbert (resp. $C^2$ Banach) manifold modeled on $H$ (resp.
$X$). Suppose that $\widetilde X\subset \widetilde H$ is dense in
$\widetilde H$, and that for each point $p\in \widetilde X$ there
exists a coordinate chart around $p$ on $\widetilde H$,
$\Phi_p:U_p\to \Phi_p(U_p)\subset H$ with $\Phi_p(p)=\theta$, such
that it restricts to a coordinate chart around $p$ on $\widetilde
H$, $\Phi_p^X: U_p\cap \widetilde X\to \Phi_p(U_p\cap \widetilde
X)\subset X$.\\
{\rm (${\bf II}$)} Let ${\cal L}:\widetilde H\to\R$ be a
continuously directional differentiable, and $F$-differentiable
functional with the following properties.
\begin{description}
\item[(a)] ${\cal L}$ satisfies the (PS) condition, and restricts to
a  $C^2$-functional on $\widetilde X$;

\item[(b)] ${\rm rank}H_k({\cal L}^b, {\cal L}^a;{\bf K})\ne 0$ for some
$k\in\N$ and regular values $a<b$;

\item[(c)] $\exists$ a finite set $\{p_1,\cdots, p_m\}\subset\widetilde X$
is contained in $K({\cal L})\cap{\cal L}^{-1}[a, b]$;

\item[(d)] Around each $p_i$ there exists a chart $\Phi_{p_i}:U_{p_i}\to \Phi_{p_i}(U_{p_i})\subset H$
as in (I) such that the functional ${\cal L}\circ
\circ(\Phi_{p_i})^{-1}: \Phi_{p_i}(U_{p_i})\to\R$ satisfies the
conditions of Theorem~2.1 of \cite{Lu1}; so $p_i$ has Morse index
$\mu_i$ and nullity $\nu_i$;

\item[(f)] Either $\mu_i>k$ or $\mu_i+ \nu_i<k$, $i=1,\cdots,m$.
\end{description}
Then ${\cal L}$ has at least one more critical point $p_0$ with
${\rm rank}C_k({\cal L}, p_0;{\bf K})\ne 0$.
\end{theorem}

\noindent{\bf Proof}. By the condition (II), the lower critical
point set of ${\cal L}$ coincides with $K({\cal L})$, and ${\cal L}$
satisfies the (PS) condition for continuous functionals. If the
conclusion is not true then $K({\cal L})=\{p_1,\cdots, p_m\}$. By
Corollary~2.6 of \cite{Lu1} and (f) we have $C_k({\cal L}, p_i)=0$,
$i=1,\cdots, m$. As in the proof of Theorem 5.1 on \cite[page
121]{Ch}  we may use  Corvellec's Morse theory for continuous
functionals \cite{Cor} to obtain a contradiction. $\Box$\vspace{2mm}

Similarly, a suitable weaker version of Theorem 5.4 on \cite[page
121]{Ch} may be given. In particular, by Theorem~2.12 of \cite{Lu1}
we have the following generalization of Corollary 5.3 therein.

\begin{theorem}\label{th:4.4}
Under the assumptions of Theorem~2.1 of \cite{Lu1}, suppose $V=H$
and the following conditions hold:
\begin{description}
\item[(i)]  ${\cal L}$ is bounded below, $F$-differentiable and
satisfies the (PS) condition;
\item[(ii)] For a small $\epsilon>0$, {\rm either} $\deg_{\rm BS}(\nabla{\cal L}, B_H(\theta,
\epsilon),\theta)=\pm 1$, {\rm or} $\deg_{\rm FPR}(A, B_X(\theta,
\epsilon),\theta)=\pm 1$ provided that the map $A$ in the condition
(F2) is $C^1$ near $\theta\in X$, where the degrees $\deg_{\rm BS}$
and $\deg_{\rm FPR}$ are as in Theorem~2.12 of \cite{Lu1}.
\end{description}
Then ${\cal L}$ has at least three critical points.
\end{theorem}

Finally,  we give a generalization of Theorem 3.12 in \cite{BaLi}.
To this goal we also need the following results, which are
Propositions 2.3 and 3.8 in \cite{BaLi}.

\begin{proposition}\label{prop:4.5}
Let a normed vector space $X$ have a direct sum decomposition
$X=X_1\oplus X_2$, where $k=\dim X_2<\infty$. For $f\in C(X,\R)$ we
have:
\begin{description}
\item[(i)]  If there exist $x_0\in X$ and $\varepsilon>0$ such that
$f(x_0+ x)>f(x_0)\;\forall x\in H_1, 0<\|x\|\le\varepsilon$, and
that $f(x_0+ x)\le f(x_0)\;\forall x\in X_2, \|x\|\le\varepsilon$,
then $C_k(f, x_0)\ne 0$. {\rm (\cite{Liu} )}

\item[(ii)]  If $f$ is bounded from below on $X_1$, and $f(x)\to -\infty$ for $x\in X_2$ as
$\|x\|\to\infty$, then $H_k(X, f^a)\ne 0$ for $a<\inf f|_{X_1}$.
{\rm (\cite[Proposition 3.8]{BaLi})}
\end{description}
\end{proposition}

\begin{theorem}\label{th:4.6}
Under the assumptions (i)-(iv) of Theorem~\ref{th:3.2}, let the
condition (i) of Proposition~\ref{prop:3.5} be satisfied. (Of course
the corresponding densely imbedded Banach spaces in $H$ are not
necessarily same, we denote by $A_\infty$ and $B_\infty$ the
corresponding maps in (i) of Theorem~\ref{th:3.2}). Let $H$ split as
$H=H^0\oplus H^+\oplus H^-$ (resp. $H=H^0_\infty\oplus
H^+_\infty\oplus H^-_\infty$) according to the spectral
decomposition of $B(x_0)$ (resp. $B_\infty(\infty)$). Let $\mu_0$,
$\nu_0$ (resp. $\mu_\infty$, $\nu_\infty$) be the Morse index and
nullity at $x_0$ (resp. infinity).
\begin{description}
\item[(I)] If (v) of Theorem~\ref{th:3.2}  and
the local linking condition as in Proposition~\ref{prop:4.5}(i) with
$X^-=H^-$ (resp. $X^-=H^0\oplus H^-$) hold, then there exists a
critical point different from $x_0$ provided $\mu_0\notin
[\mu_\infty, \mu_\infty+ \nu_\infty]$ (resp. $\mu_0+ \nu_0\notin
[\mu_\infty, \mu_\infty+\nu_\infty]$).

\item[(II)] If (ii)-(iv) and (${\bf AC^+}$) of  Proposition~\ref{prop:3.5} hold,
and ${\cal L}$ is bounded from below on $H^0_\infty\oplus
H^+_\infty$, then there exists a nontrivial critical point provided
$\mu_\infty\ne\mu_0$.
\end{description}
\end{theorem}

\noindent{\bf Proof}. (I) Otherwise we have $K({\cal L})=\{x_0\}$.
By Proposition~\ref{prop:3.3}(ii), $C_\ast({\cal L},\infty)\cong
C_\ast({\cal L}, x_0)$. Moreover, Proposition~\ref{prop:4.5}(i)
implies $C_{\mu_0}({\cal L}, x_0)\ne 0$. Hence $C_{\mu_0}({\cal
L},\infty)\ne 0$. This gives a contradiction by
Theorem~\ref{th:3.2}.

(II) Similarly, assume $K({\cal L})=\{x_0\}$. By
Proposition~\ref{prop:3.3}(ii), $C_\ast({\cal L},\infty)\cong
C_\ast({\cal L}, x_0)$. Proposition~\ref{prop:3.5}(a) yields
$C_k({\cal L}, x_0)=\delta_{k\mu_0}{\bf K}$.  By
Proposition~\ref{prop:4.5}(ii), $C_{\mu_\infty}({\cal L},\infty)\ne
0$. This contradiction proves the desired conclusion. $\Box$
\vspace{2mm}

.

\section{Critical groups of sign-changing critical points}\label{sec:5}
\setcounter{equation}{0}

In this section we shall present the corresponding version of the
results on critical groups of sign-changing critical points in
\cite{BaChWa} and \cite{LiuWaWe} in our framework and sketch how to
prove them with our results in the previous sections. It is also
possible to generalize some of \cite{Ba1}. They are left to the
interested reader.

Let the Hilbert space $(H, (\cdot,\cdot)_H)$ and the Banach space
$(X, \|\cdot\|_X)$ satisfy the condition ({\bf S}) as in \cite[\S
2.1]{Lu1}. Let $P_H$ be a closed convex in $H$ so that
$P\stackrel{\triangle}{=}P_H\cap X$, a closed convex cone in $X$,
satisfies:
\begin{equation}\label{e:5.1}
\hbox{(i) ${\rm int}_X(P)\ne\emptyset$, (ii) $\exists\;e\in {\rm
int}_X(P)$ with $(u,e)_H>0$ for all $u\in P\setminus\{\theta\}$.}
\end{equation}
Then $H$ (resp. $X$) is  partially ordered by by $P_H$ (resp. $P$).
For $u,v\in H$ (resp. $X$) we write: $u\ge v$ if $u-v\in P_H$ (resp.
$u-v\in P$); $u>v$ if $u-v\in P_H\setminus\{\theta\}$ (resp. $u-v\in
P\setminus\{\theta\}$). When $u, v\in X$ we also write $u\gg v$ if
$u-v\in {\rm int}_X(P)$. A map $f:H\to H$ (resp. $f:X\to X$) is
called  order preserving if $u\ge v\Rightarrow f(u)\ge f(v)\;\forall
u, v\in H$ (resp. $X$). In particular, $f:X\to X$ is said to be
strongly order preserving if $u> v\Rightarrow f(u)\gg f(v)\;\forall
u, v\in X$.

Recall that  above Proposition~\ref{prop:2.1} we have showed that a
continuously directional differentiable map is  locally Lipschitz
continuous  and strictly G-differentiable.

Let the assumptions of Theorem~2.1 of \cite{Lu1} hold for $V=H$. We
also assume:

\begin{description}
\item[(${\bf L}_0$)] The assumptions of Proposition~2.26 of
\cite{Lu1} hold for $V=H$, that is, the assumptions of Theorem~2.1
of \cite{Lu1} hold for $V=H$,  the map $A: X\to X$ in the condition
({\bf F2}) of \cite{Lu1} is Fr\'echet differentiable, and there
exist positive constants $\eta_0'$ and $C'_2>C'_1$ such that
$$
C'_2\|u\|^2\ge (P(x)u, u)\ge C'_1\|u\|^2\quad\forall u\in
H,\;\forall x\in B_H(\theta,\eta_0')\cap X.
$$

\item[(${\bf L}_1$)] ${\cal L}:H\to\R$ is $C^{2-0}$, and satisfies the
$(PS)_c$ condition for any $c\in\R$. Moreover, all critical points
of ${\cal L}$ are contained in $X$, and $H$ and $X$ induce an
equivalent topology on the critical set of ${\cal L}$ \footnote{The
final assumption is used in the proof of Claim 3, cf. the proof of
Theorem~\ref{th:5.1}. There is no such a assumption in
\cite{BaChWa}.}

\item[(${\bf L}_2$)] The map  $id_X-A:X\to X$ is strongly order preserving.

\item[(${\bf L}_3$)] The smallest eigenvalue of $B(\theta)$, which
is equal to $\inf_{\|u\|=1}(B(\theta)u,u)_H$ by \cite[Prop.6.9]{Bre}
and Proposition~B.2 of \cite{Lu1}, is simple and its eigenspace
(contained in $X$ by ({\bf D1})) is spanned by a positive
eigenvector (i.e. sitting in ${\rm Int}_X(P)$).

\item[(${\bf L}_4$)] One of the following holds:
\begin{description}
\item[(i)] ${\cal L}$ is bounded below;

\item[(ii)] For every $u\in H\setminus\{\theta\}$ it holds that ${\cal
L}(tu)\to-\infty$ as $t\to\infty$. Moreover, there exists $\flat<0$
such that ${\cal L}(u)\le\flat$ implies $D{\cal L}(u)(u)<0$;

\item[(iii)] There exist a compact self-adjoint linear operator
$Q_\infty\in L(H)$ and a positive definite operator $P_\infty\in
L(H)$ such that $\nabla{\cal L}(u)=P_\infty u-Q_\infty u+ o(\|u\|)$
as $\|u\|\to\infty$, where $\nabla{\cal L}$ is the gradient of
${\cal L}$. The smallest eigenvalue of $B_\infty:=P_\infty-Q_\infty$
is simple and its eigenspace is spanned by a positive eigenvector
$e_\infty\in{\rm Int}_X(P)$ such that $(u, e_\infty)_H>0$ for every
$u\in P\setminus\{\theta\}$. Moreover, $B_\infty|_X\in L(X)$, and if
a subset $S\subset X$ is bounded in $H$ then $A(S)$ is also bounded
in $H$.
\end{description}
\end{description}

Let $\mu_0:=\dim H^-$ and $\nu_0:=\dim H^0$. Define
$\mu_\infty=\nu_\infty=0$ in the case (i) of (${\bf L_4}$), and
$\mu_\infty=\infty$ and $\nu_\infty=0$ in the case (ii) of (${\bf
L_4}$). For the case (iii) of (${\bf L_4}$), let $\mu_\infty$ be the
number of negative eigenvalues of $B_\infty$ (counted with
multiplicities) and $\nu_\infty=\dim {\rm Ker}(B_\infty)$. An
element $x\in X$ is called a {\bf subcritical} (resp. {\bf
supercritical}) critical point of ${\cal L}$ if $\nabla{\cal
L}(x)\le 0$ (resp. $\nabla{\cal L}(x)\ge 0$).

By ({\bf F2}) the map $A:X\to X$ is continuously directional
differentiable. It follows from Proposition~B.1(ii) of \cite{Lu1}
that $A$ is a locally Lipschitz map from $X$ to $X$. Note that
$\nabla{\cal L}$ is $C^{1-0}$, and equal to $A$ on $X$ by (${\bf
L}_2$), i.e. $\nabla{\cal L}(x)=A(x)\;\forall x\in X$.

Consider the negative gradient flow $\varphi(t, x)$ of  ${\cal L}$
on $H$ defined by
\begin{equation}\label{e:5.2}
\frac{d}{dt}\varphi(t, x)=-\nabla{\cal
L}(\varphi(t,x))\quad\hbox{and}\quad \varphi(0,x)=x\in H.
\end{equation}
Let $\varphi^t(\cdot)=\varphi(t,\cdot)$. As in \cite{BaChWa} it
restricts to a continuous local flow on $X$, still denoted by
$\varphi^t(\cdot)$, and $t\mapsto {\cal L}(\varphi^t(x))$ is
strictly decreasing for any $x\notin K({\cal L})$. Let ${\cal
D}=P\cup(-P)$ and ${\cal D}^\ast={\cal D}\setminus\{\theta\}$, and
let ${\cal L}^d=\{{\cal L}(x)\le d\}$ for $d\in\R$. The following
result is a generalization of Theorem 3.4 in \cite{BaChWa}.

\begin{theorem}\label{th:5.1}
Under the above assumptions (${\bf L_0}$)-(${\bf L_4}$), if
$\mu_0\ge 2$ and $\mu_\infty+ \nu_\infty\le 1$ then ${\cal L}$ has a
sign-changing critical point $x_1$ with ${\cal L}(x_1)<0$. If all
sign changing critical points (i.e. those in $X\setminus{\cal D}$ )
with critical values contained in a bounded interval of $(-\infty,
0)$ are isolated, then there exists a sign changing critical point
$x_1$ with ${\cal L}(x_1)<0$ and $C_1({\cal L}, x_1;{\bf K})\ne 0$.
Furthermore, if near $x_1$ the conditions of Theorem~2.1 of
\cite{Lu1} hold and $x_1$ has nullity $1$ then $x_1$ is of mountain
type and $C_k({\cal L}, x_1;{\bf K})\cong\delta_{k1}{\bf F}$. {\rm
(Note: It is sufficient that ${\cal L}$ satisfies the $(PS)_c$
condition for each $c<0$).}
\end{theorem}

\noindent{\bf Proof}. Firstly, we prove:

\noindent{\it Claim 1}. \textsf{Suppose that all sign changing
critical points with critical values contained in a bounded interval
of $(-\infty, 0)$ are isolated. Then for any interval $[a, b]\subset
(-\infty, 0)$ there exist only finitely many sign changing critical
points with critical values in $[a,b]$.}

Otherwise, let $\{u_n\}$ be infinite such points with $\{{\cal
L}(u_n)\}\subset [a, b]\subset (-\infty, 0)$. We may assume ${\cal
L}(u_n)\to c$ (by passing to a subsequence if necessary). By the
$(PS)_c$ condition we may assume $u_n\to u_0$ in $H$, and $u_n\to
u_0$ in $X$ \textbf{because of the final assumption in (${\bf
L_1}$)}. Clearly, $u_0\ne \theta$ because of ${\cal L}(u_0)\le b<0$.
Since each $u_n$ sits in an open subset $H\setminus(P_H\cup(-P_H))$,
it belongs to
$X\cap(H\setminus(P_H\cup(-P_H)))=X\setminus(P\cup(-P))$ as well.
Then either $u_0\in X\setminus(P\cup(-P))$ or $u_0\in\partial {\cal
D}\setminus\{\theta\}$. In the first case $u_0$ also sits in
$H\setminus(P_H\cup(-P_H))$. This contradicts the assumption that
all sign changing critical points with negative critical values are
isolated. In the latter case either $u_0\in\partial
P\setminus\{\theta\}$ or $u_0\in (-\partial P)\setminus\{\theta\}$.
Note that $id_X-A$ is strongly order preserving. If $u_0\in
\partial P\setminus\{\theta\}$ then $u_0=u_0-A(u_0)\gg \theta- A(\theta)=\theta$, i.e.
$u_0\in{\rm Int}_X(P)$. This leads to a contradiction. If $u_0\in
\partial (-P)\setminus\{\theta\}$, i.e. $-u_0\in
\partial P\setminus\{\theta\}$ then
$$
\theta=u_0+ (-u_0)- A(u_0+ (-u_0))\gg u_0-A(u_0)=u_0
$$
because $u_0+ (-u_0)\ge u_0$. Hence $-u_0\in{\rm Int}_X(P)$, which
yields a contradiction again. Claim 1 is proved.

Next, following the methods of proof of \cite[Th.3.4]{BaChWa}, since
$id_X-A:X\to X$ is strongly order preserving by (${\bf L}_2$) we get
that $\varphi^t(x)\in{\rm Int}_X({\cal D})$ for all $x\in{\cal
D}^\ast$ and $t>0$. For $d\in\R$ let $(j_d)_1$ be the homomorphism
from $H_1(X, {\cal D}^\ast;{\bf K})$ to $H_1(X, {\cal L}^d\cup{\cal
D}^\ast;{\bf K})$ induced by the inclusion $j_d: (X, {\cal
D}^\ast)\hookrightarrow (X, {\cal L}^d\cup{\cal D}^\ast)$. Let
$$
\Gamma:=\{d\in\R\,|\, (j_d)_1\ne 0\}\quad\hbox{and}\quad
c:=\sup\Gamma.
$$

Since $\mu_\infty+ \nu_\infty\le 1$, either (i) of (${\bf L}_4$) or
(iii) of (${\bf L}_4$) occurs. Let $v_\infty=e/\|e\|$ in case (i),
and $v_\infty=e_\infty/\|e_\infty\|$ in case (iii). We conclude:

\noindent{\it Claim 2}.  \textsf{${\cal L}$ is bounded below on
$H_1:=H\cap\langle v_\infty\rangle^\bot$.}

We only need to prove this in the latter case. Since the smallest
eigenvalue of $B_\infty$, given by $\inf_{\|u\|=1}(B_\infty u,
u)_H=(B_\infty v_\infty, v_\infty)_H$, is simple, we have
$$
\lambda_2:=\inf\{(B_\infty u, u)_H\,|\, \|u\|=1,\;u\in
H_1\}>\lambda_1:=(B_\infty v_\infty, v_\infty)_H.
$$
Obverse that $e_\infty\in X$ implies that $X_1:=X\cap\langle
v_\infty\rangle^\bot$ is dense in $H_1$.

Since $\mu_\infty+\nu_\infty\le 1$ we have three cases: (a)
$\mu_\infty=0=\nu_\infty$, (b) $\mu_\infty=0$ and $\nu_\infty=1$,
(c) $\mu_\infty=1$ and $\nu_\infty=0$. They corresponds to
$\lambda_1>0$, $\lambda_1=0$,  and $\lambda_1<0$ but $\lambda_2>0$,
respectively. So we always have $\lambda_2>0$.

Take a large $N>0$ so that  $\|\nabla{\cal L}(u)-B_\infty u\|<
\frac{\lambda_2}{4}\|u\|$ as $\|u\|\ge N$. For any $u\in H_1$ with
$\|u\|>N$ let $\bar u=N\cdot u/\|u\|$. By the assumption
$A\bigl(X\cap B_H(\theta, N)\bigr)$ is bounded in $H$. So there
exists a $M>0$ such that for any $u\in X\cap B_H(\theta, N)$,
$$
|{\cal L}(u)|=|{\cal L}(u)-{\cal L}(\theta)|=\left|\int^1_0D{\cal
L}(tu)(u)dt\right|=\left|\int^1_0(A(tu), u)_Hdt\right|\le MN.
$$
This also holds for all $u\in  B_H(\theta, N)$ because $X\cap
B_H(\theta, N)$ is dense in $B_H(\theta, N)$. On the other hand
\begin{eqnarray*}
{\cal L}(u)-{\cal L}(\bar u)&=&\int^1_0(\nabla{\cal L}(t u+
(1-t)\bar u), (u-\bar
u))_Hdt\\
&\ge& \int^1_0(B_\infty(t u+ (1-t)\bar u), u-\bar
u)_Hdt-\frac{\lambda_2}{4}\|u-\bar u\|\cdot\|t u+ (1-t)\bar
u\|\\
&\ge &\frac{1}{2}(B_\infty(u-\bar u), u-\bar u)_H+(B_\infty(\bar u),
u-\bar u)_H -\frac{\lambda_2}{4}\|u-\bar u\|\cdot\|\bar
u\|\\
&\ge &\frac{\lambda_2}{2}\|u-\bar u\|^2+(B_\infty(\bar u), u-\bar
u)_H -\frac{\lambda_2}{4}\|u-\bar u\|\cdot\|\bar u\|.
\end{eqnarray*}
Claim 2 follows immediately.

This implies that any $d<\inf{\cal L}|_{X_1}=\inf{\cal L}|_{H_1}$
belongs to $\Gamma$. Hence $c$ is finite. Since $\mu_0\ge 2$ the two
smallest eigenvalues of $B(\theta)$ are negatives and the
corresponding eigenspaces are contained in $X$ by ({\bf D1}). Let
$e_1$ and $e_2$ two normalized eigenvectors belonging to the two
smallest eigenvalues of $B(\theta)$. By (${\bf L}_3$), $e_1\in{\rm
Int}_X(P)\subset{\rm Int}_X({\cal D})$. Let $S_\rho$ be the sphere
of radius $\rho$ in ${\rm Span}\{e_1, e_2\}$. It easily follows from
Proposition~2.26 of \cite{Lu1} that $\max{\cal L}(S_\rho)<0$.
\footnote{This is only place where Proposition~2.26 of \cite{Lu1} is
used. The assumptions of Theorem~2.1 of \cite{Lu1} and (${\bf
L_1}$)-(${\bf L_4}$) are sufficient to other arguments.} As in the
proof of \cite[Lem.4.2]{BaChWa} we can prove $c<0$.

\noindent{\it Claim 3}. \textsf{ $c$ is  a critical value of ${\cal
L}$, and hence ${\cal L}$ has a sign-changing critical point with
negative critical value.}

This may be proved by a standard deformation argument. In view of
Claim 1 let us suppose that there exist only finitely many sign
changing critical points $x_1,\cdots, x_q$ at the level $c$. Note
that Claim 1 also implies that there exist a $\eta>0$ such that no
number in $[c-\eta, c+ \eta]\setminus\{c\}$ is a critical value of
sign changing critical points. Repeating the remainder of proof of
\cite[Th.3.4]{BaChWa} we get some $i\in\{1,\cdots,q\}$ such that
$C_1({\cal L}, x_i;{\bf K})\ne 0$. The final conclusions follow from
Corollary~2.9(ii) of \cite{Lu1}. $\Box$

Corresponding with Theorem 3.5 in \cite{BaChWa} we have:

\begin{theorem}\label{th:5.2}
Under the assumptions (${\bf L_0}$)-(${\bf L_4}$) above, if
$\mu_0\ge 2$ and there exist a subcritical critical point
$\underline{x}$ and a supercritical critical point $\bar x$ of
${\cal L}$ such that $\underline{x}\ll\theta\ll\bar x$, then the
conclusions of Theorem~\ref{th:5.1} is still true.
\end{theorem}

Similarly, we can get the corresponding results with Theorems 3.6,
3.8 in \cite{BaChWa} as follows.

\begin{theorem}\label{th:5.3}
Under the assumptions  (${\bf L_0}$)-(${\bf L_4}$) above, if
$\mu_\infty\ge 2$ and $\mu_0+ \nu_0\le 1$ then ${\cal L}$ has a
sign-changing critical point $x_1$ with ${\cal L}(x_1)>0$. If all
sign changing critical points (i.e. those in $X\setminus{\cal D}$ )
with critical values contained in a bounded interval of $[0,
\infty)$ are isolated, then there exists a sign changing critical
point $x_1$ with ${\cal L}(x_1)>0$, Morse index $\mu\in\{1,2\}$ and
$C_0({\cal L}, x_1;{\bf K})=0=C_1({\cal L}, x_1;{\bf K})$,
$C_2({\cal L}, x_1;{\bf K})\ne 0$.  Furthermore, if near $x_1$ the
conditions of Theorem~2.1 of \cite{Lu1} hold then $x_1$ is neither a
local minimum nor of mountain pass type,  and $C_k({\cal L},
x_1;{\bf K})=\delta_{k2}{\bf K}$ holds for all $k$ provided that
$\mu=2$ or the nullity $\nu\le 1$.
\end{theorem}

\begin{theorem}\label{th:5.4}
Under the assumptions  (${\bf L_0}$)-(${\bf L_4}$) above, let
$\nu_0=0=\nu_\infty$, $\mu_\infty\ge 1$ and $\mu_0\ne \mu_\infty$.
Suppose also that  all sign changing critical points  are isolated.
Then
\begin{description}
\item[(i)] If $\mu_0\ge 1$ then ${\cal L}$ has  a sign changing critical point $x_1$ which satisfies
either ${\cal L}(x_1)>0$ and $C_{\mu_0+1}({\cal L}, x_1;{\bf K})\ne
0$ or ${\cal L}(x_1)<0$ and $C_{\mu_0-1}({\cal L}, x_1;{\bf K})\ne
0$.

\item[(ii)] If (${\bf L_4}$) (iii) applies with $\mu_\infty\ge 2$
then ${\cal L}$ has  a sign changing critical point $x_1$ with
$C_{\mu_\infty}({\cal L}, x_1;{\bf K})\ne 0$.
\end{description}
\end{theorem}

The last theorem can be obtained by completely repeating the proof
of Theorem 3.8 in \cite{BaChWa}.

\noindent{\bf Proof of Theorem~\ref{th:5.3}}. For reader's
convenience we follow the proof ideas of \cite[Theorem 3.6]{BaChWa}
to give necessary details.
 Since $\mu_\infty\ge 2$ either (${\bf L_4}$)-(ii) or (${\bf
L_4}$)-(iii) occurs.

\noindent{\it Claim 4}. \textsf{There exist two orthogonal unit
vectors $v_\infty\in{\rm Int}_X(P)$ and $u_\infty\in X$ such that
${\cal L}(u)<0$ for $u\in {\rm span}\{v_\infty, u_\infty\}$ with
$\|u\|\ge R$.}

In fact, In the latter case, the two smallest eigenvalue of
$B_\infty$,
$$
\lambda_1:=(B_\infty v_\infty,
v_\infty)_H<\lambda_2:=\inf\{(B_\infty u, u)_H\,|\, \|u\|=1,\;u\in
H_1\}
$$
are negative, where $v_\infty=e_\infty/\|e_\infty\|$ and
$H_1:=H\cap\langle v_\infty\rangle^\bot$. Since $X$ is dense in $H$,
$X\cap H_1\ne\emptyset$. Let $u_\infty$ be a unit vector in $X\cap
H_1$. Then $(v_\infty, u_\infty)_H=0$.

As in the arguments below Claim 2  in the proof of
Theorem~\ref{th:5.1} we have $N>0$ and $M>0$ such that
\begin{eqnarray*}
&&\|\nabla{\cal L}(u)-B_\infty u\|<
\frac{|\lambda_2|}{4}\|u\|\quad\hbox{as}\;\|u\|\ge N,\\
&& |{\cal L}(u)|\le MN\quad\forall u\in B_H(\theta, N).
\end{eqnarray*}
For any $u\in {\rm span}\{v_\infty, u_\infty\}$ with $\|u\|>N$ let
$\bar u=N\cdot u/\|u\|$. Then
\begin{eqnarray*}
{\cal L}(u)&=&\int^1_0(\nabla{\cal L}(t u+ (1-t)\bar u), (u-\bar
u))_Hdt+ {\cal L}(\bar u)\\
&\le& \int^1_0(B_\infty(t u+ (1-t)\bar u), u-\bar
u)_Hdt+\frac{|\lambda_2|}{4}\|u-\bar u\|\cdot\|t u+ (1-t)\bar
u\|+ MN\\
&\le &\frac{1}{2}(B_\infty(u-\bar u), u-\bar u)_H+(B_\infty(\bar u),
u-\bar u)_H +\frac{|\lambda_2|}{4}\|u-\bar u\|\cdot\|\bar
u\|+ MN\\
&\le &\frac{\lambda_2}{2}\|u-\bar u\|^2+(B_\infty(\bar u), u-\bar
u)_H +\frac{|\lambda_2|}{4}\|u-\bar u\|\cdot\|\bar u\|+ MN.
\end{eqnarray*}
So Claim 4 follows from this in this case.

In the former case take any  unit vector $v_\infty\in{\rm
Int}_X(P)$.  As above we can choose another unit vector $u_\infty\in
X$ which is  orthogonal to $v_\infty$.  Since for any $u\in H$ with
${\cal L}(u)\le a$ it holds that $D{\cal L}(u)(u)<0$, ${\cal L}^a$
is a manifold with $C^1$-boundary ${\cal L}^{-1}(a)$, and ${\cal
L}^{-1}(a)$ is transversal to the radial vector field. Moreover, for
any $u\in H\setminus\{\theta\}$, ${\cal L}(tu)\to-\infty$ as
$t\to\infty$. So for each $u\in \partial B_H(\theta, 1)$ there
exists a $t_u\in (0, \infty)$ such that ${\cal L}(t_u\cdot u)=a$ and
that $u\to t_u$ is continuous by the implicit function theorem. This
shows that the map $\partial B_H(\theta, 1)\to {\cal
L}^{-1}(a),\;u\to t_u\cdot u$ is a homeomorphism. It follows that
$\partial B_H(\theta, 1)\cap {\rm span}\{v_\infty, u_\infty\}$ is
compact and hence  $\partial B_H(\theta, 1)\cap {\rm
span}\{v_\infty, u_\infty\}\subset B_H(\theta, R)$ for some $R>0$.
We have still  ${\cal L}(u)<0$ for $u\in {\rm span}\{v_\infty,
u_\infty\}$ with $\|u\|\ge R$. That is, Claim 4 holds.

For $R$ in Claim 4 let us set
\begin{eqnarray*}
&&B_R:=\{sv_\infty+ tu_\infty:\;|s|\le R,\;0\le t\le R\},\\
&&\partial B_R:=\{sv_\infty+ tu_\infty:\;|s|=R,\;\hbox{or}\;
t\in\{0,1\}\}.
\end{eqnarray*}
Then $\partial B_R\subset{\cal L}^0\cup{\cal D}$. Let
$\beta:=\max{\cal L}(B_R)$ and $\xi_\beta\in H_2({\cal
L}^\beta\cup{\cal D}, {\cal L}^0\cup{\cal D};{\bf K})$ be the image
of $1\in{\bf K}\cong H_2(B_R, \partial B_R;{\bf K})$ under the
homomorphism
$$
{\bf K}\cong H_2(B_R, \partial B_R; {\bf K})\to H_2({\cal
L}^\beta\cup{\cal D}, {\cal L}^0\cup{\cal D};{\bf K})
$$
induced by the inclusion $(B_R, \partial B_R)\hookrightarrow B_R,
\partial B_R)$. For $\gamma\le\beta$ let
$$
(j_\gamma)_2: H_2({\cal L}^\gamma\cup{\cal D}, {\cal L}^0\cup{\cal
D};{\bf K})\to H_2({\cal L}^\beta\cup{\cal D}, {\cal L}^0\cup{\cal
D};{\bf K})
$$
be the homomorphism induced by the inclusion. Set
\begin{equation}\label{e:5.3}
\Gamma:=\{\gamma\le\beta\,|\, \xi_\beta\in{\rm
Image}(j_\gamma)_2)\}\quad\hbox{and}\quad c:=\inf\Gamma.
\end{equation}

Let $e_1\in{\rm Int}_X(P)$ be the eigenvector of $B(\theta)$
belonging to the first eigenvalue
$\lambda_1=\inf_{\|u\|=1}(B(\theta)u,u)_H$, and let $X_1=\langle
e_1\rangle$ and $X_2:=X_1^\bot\cap X$. Since $\mu_0+ \nu_0\le 1$.
There are three cases: (a) $\mu_0=0=\nu_0$, (b) $\mu_0=1$ and
$\nu_0=0$, (c) $\mu_0=0$ and $\nu_0=1$. For the first two case we
have $X_2\subset H^+$. In the third case, $\lambda_1=0$ and
$X_1={\rm Ker}(B(\theta))$, and $H^-=\{\theta\}$. We also get
$X_2\subset H^+$. It follows from \cite[(2.74)]{Lu1} that
$$
{\cal L}(u)\ge \frac{a_1}{4}\|u\|^2
$$
for all $u\in B_H(\theta, \rho_0)\cap H^+$. Hence for some small
$\rho>0$ it always holds that $\inf\{{\cal L}(u)\,|\, u\in
X_2,\,\|u\|=\rho\}>0$. This is what is needed in the proof of
\cite[Lem.4.3]{BaChWa}. It leads to $\xi_\beta\ne 0$ and hence
$\beta\in\Gamma$. Moreover,  $(j_0)_2=0$ implies that
$0\notin\Gamma$.

Since we have assumed that $\theta$ is an isolated critical
point\footnote{In \cite{BaChWa} it was claimed that since
$\mu_0+\nu_0\le 1$ the sign changing solutions cannot accumulate at
$0$.}, and that all sign changing critical points with critical
values in a bounded interval of $[0, \infty)$, by the proof of Claim
3 in the proof of Theorem~\ref{th:5.1} we can derive that there
exist only finitely many sign changing critical points with critical
values in $[0, \beta]$. It follows that ${\cal L}^0\cup{\cal D}$ is
a strong deformation retract of ${\cal L}^\gamma\cup{\cal D}$ for
$\gamma>0$ small enough. So $c>0$.

The remained arguments are the same as those of
\cite[Th.3.6]{BaChWa} (as long as slightly modifications as in the
proof of Theorem~\ref{th:5.1}. $\Box$

\begin{remark}\label{rm:5.5}
{\rm Let us outline a possible way to weaken the conditions of
Theorem~\ref{th:5.1}, that is, removing the assumption that
${\cal L}$ is $C^{2-0}$, but adding the condition\\
$({\bf 6.5.1})$ ``For any  subset $S\subset X$, which is bounded in
$H$,
the image $A(S)$ is bounded in $X$'' in case (${\bf L_4}$)-(i); \\
$({\bf 6.5.2})$ ``For any $c\in\R$ and small $\varepsilon>0$,
$A(X\cap{\cal L}^{-1}([c-\varepsilon, c+\varepsilon]))$ is bounded
in $X$'' in case (${\bf L_4}$)-(iii).

Since ({\bf F2}) and Proposition~B.1(ii) of \cite{Lu1} imply that
the map $A:X\to X$ is locally Lipschitz, we get a (local) flow on
$X\setminus K({\cal L})$,
\begin{equation}\label{e:5.4}
\left.\begin{array}{ll}
 \frac{d}{dt}\sigma(t,
x)=- A(\sigma(t,x))\\
\sigma(0,x)=x\in X\setminus K({\cal L}),
\end{array}\right\}
\end{equation}
where $K({\cal L})$ is the critical set of ${\cal L}$. By (${\bf
L}_2$) we get that $\sigma(t,x)\in{\rm Int}_X({\cal D})$ for all
$x\in{\cal D}^\ast$ and $t>0$. The key is how to prove Claim 3 in
the present assumptions. Note that we have proved $c<0$ above Claim
3. Then Claim 1 implies that ${\cal L}^{-1}(c)$ contains at most
finitely many sign changing critical points $x_1,\cdots, x_q$ and
 that there exist a
$\eta>0$ such that no number in $[c-\eta, c+ \eta]\setminus\{c\}$ is
a critical value of sign changing critical points.

 By contradiction, suppose that $c$ is not a critical value of
${\cal L}$.  Then  the PS condition implies that there exist
$\varepsilon>0$ and $\delta>0$ such that
$$
\|\nabla{\cal L}(u)\|=\|D{\cal L}(u)\|\ge\delta\quad\forall
u\in{\cal L}^{-1}[c-\varepsilon, c+ \varepsilon].
$$
We may assume $\varepsilon<\eta$. By (${\bf L_2}$) and
$\|u\|_X\ge\|u\|\;\forall u\in X$ we get
\begin{equation}\label{e:5.5}
\|A(u)\|_X\ge\|A(u)\|=\|D{\cal L}(u)\|\ge\delta\quad\forall u\in
X\cap {\cal L}^{-1}[c-\varepsilon, c+ \varepsilon].
\end{equation}
For $x\in X\cap {\cal L}^{-1}[c-\varepsilon, c+ \varepsilon]$ let
$[0, T_x)$ be the maximal existence interval of the flow in
(\ref{e:5.4}) on $X\cap {\cal L}^{-1}[c-\varepsilon, c+
\varepsilon]$. Then
$$
\eta+\varepsilon\ge{\cal L}(x)-{\cal
L}(\sigma(t,x))=\int^t_0\|A(\sigma(s,x))\|^2ds\ge \delta^2t
$$
for any $t\in [0, T_x)$. So $T_x\le (\eta+\varepsilon)/\delta^2$.

{\bf In case (${\bf L_4}$)-(i)}, ${\cal L}$ is coercive by a result
of \cite{CaLiWi}. It follows that ${\cal L}^{c+\varepsilon}$ is
bounded in $H$ and hence $A(X\cap{\cal L}^{c+\varepsilon})$ is
bounded in $X$ by the assumption (6.5.1). Namely, there exists a
$N>0$ such that $\|A(x)\|_X\le N$ for all $x\in X\cap{\cal
L}^{c+\varepsilon}$. Then
$$
{\rm dist}_X(\sigma(t_2,x),
\sigma(t_1,x))\le\int^{t_2}_{t_1}\left\|\frac{d}{dt}\sigma(t,
x)\right\|_Xdt=\int^{t_2}_{t_1}\|A(\sigma(t,x))\|_X dt\le N(t_2-t_1)
$$
for any $0\le t_1<t_2<T_x$. It follows that the limit $\lim_{t\to
T_x-0}\sigma(t,x)$ exists in $X$ and
$$
{\cal L}\left(\lim_{t\to T_x-}\sigma(t,x)\right)=c-\varepsilon.
$$
As usual we can use $\sigma$ to construct a deformation retract from
${\cal L}_X^{c+\varepsilon}\cup {\cal D}^\ast$ to ${\cal
L}_X^{c-\varepsilon}\cup {\cal D}^\ast$, where ${\cal
L}_X^d:=X\cap{\cal L}^d$ for $d\in\R$. This is a contradiction.
Hence $c$ is a critical value.

{\bf In case (${\bf L_4}$)-(iii)}, by the assumption (6.5.2),
$A(X\cap{\cal L}^{-1}([c-\varepsilon, c+\varepsilon]))$ is bounded
in $X$ for some small $\varepsilon>0$. Then the same method leads to
a contradiction yet.}$\Box$
\end{remark}

 Finally, we are going to generalize the following result, which is an abstract summary
 of the arguments in \cite{LiuWaWe}.

\begin{theorem}\label{th:5.6}
Let $H$ be a Hilbert space, and let  $P_H\ne H$ be a closed cone,
i.e. $P_H=\bar P_H$ is convex and satisfies $\R^+\cdot P_H\subset
P_H$, $P_H\cap(-P_H)=\{\theta\}$. Suppose that  a $C^2$-functional
   ${\cal L}:H\to\R$ has critical point $\theta$ and satisfies the following properties.
\begin{description}
\item[(i)] ${\cal L}$ is bounded from below, and satisfies the (PS)
condition. \footnote{By a result of \cite{CaLiWi} these two
conditions imply the coercivity of ${\cal L}$. }

\item[(ii)] $P_H$ is positively invariant under the negative gradient
flow $\varphi^t$ of ${\cal L}$,
$$
 \frac{d}{dt}\varphi^t(u)=- \nabla{\cal L}(\varphi^t(u))\quad\hbox{and}\quad \varphi^0(u)=u\in H.
$$

\item[(iii)] There exists a positive element $e\in P_H\setminus\{\theta\}$ such that the cone
$$
\hbox{$D:=\{u\in H\,|\, u\ge e\}\subset P_H$
 (resp. $-D$)}
$$
 contains all positive (resp. negative) critical
points of $\Phi$. (Note that $D\cap(-D)=\emptyset$). Let
$D_\varepsilon:=\{u\in H\,|\, {\rm dist}(u,D)\le\varepsilon\}$ for
$\varepsilon>0$.

\item[(iv)] There exists a $\varepsilon_0>0$ such that
for each $\varepsilon\in (0, \varepsilon_0]$,
$D_{\varepsilon_0}\cap(-D_{\varepsilon_0})=\emptyset$\footnote{By
considering the functional $f(x)=(x,e)_H$ we obtain $D\subset\{f>
0\}$ and $-D\subset\{f<0\}$. This implies ${\rm dist}(D, -D)\ge
2\|e\|$ and thus ${\rm dist}(D_{\varepsilon_0},
-D_{\varepsilon_0})\ge 2\|e\|-2\varepsilon_0$ for
$\varepsilon<\|e\|$!}, and
 $D_\varepsilon$ and $-D_\varepsilon$  are strictly positively invariant for the the
 flow $\varphi^t$.
\end{description}
Let $W_\varepsilon:=D_\varepsilon\cup(-D_{\varepsilon})$ and let
$i_c:({\cal L}^c\cup W_\varepsilon, W_\varepsilon)\to (H,
W_\varepsilon)$ be the inclusion. Then
\begin{eqnarray}\label{e:5.6}
&&c_1:=\inf\{c\in\R\,|\, i^\ast_c:\bar H^{1}(H,
W_\varepsilon;\Z_2)\to
\bar H^{1}({\cal L}^c\cup W_\varepsilon, W_\varepsilon;\Z_2)\nonumber\\
&&\hspace{70mm}\hbox{is a monomorphism}\}
\end{eqnarray}
is finite and $K^\ast_{c_1}:=\{u\in H\setminus W_\varepsilon,|\,
{\cal L}(u) = c_1,\;{\cal L}'(u)=0\}\ne\emptyset$ for
$0<\varepsilon\le\varepsilon_0$, where $i^\ast_c$ is induced by the
inclusion $i_c$. Moreover, ${\cal L}$ has at least two nontrivial
critical points in $H\setminus W_\varepsilon$ provided that
$c_1<{\cal L}(\theta)$ and
\begin{eqnarray}
&&\hbox{each $u\in K^\ast_{c_1}$ with Morse index $\mu(u)=0$ has nullity $\nu(u)\le 1$},\label{e:5.7}\\
&&\hbox{$C_q({\cal L}, \theta;\Z_2)\cong\delta_{qn}\Z_2$ for some
$n\ge 2$ and any $q\in\Z$}.\label{e:5.8}
\end{eqnarray}
\end{theorem}

\begin{remark}\label{rm:5.7}
{\rm (i) Actually, the assumptions (i)-(ii) in Theorem~\ref{th:5.6}
imply that ${\cal L}$ has a critical point $u\in P_H$ with ${\cal
L}(u)=\inf_{v\in P_H}{\cal L}(v)$ (\cite[Theorem 2.1]{LiuSun}). In
the same way the assumption (iv) yields a critical point sitting in
$-D\subset P_H$.\\
(ii) Note that (\ref{e:5.8}) holds if $\theta$ is a nondegenerate
critical point of ${\cal L}$ with Morse index $n\ge 2$.\\
(iii) If  ${\cal K}:=id-\nabla{\cal L}: H\to H$ maps $D_\varepsilon$
(resp. $-D_\varepsilon$) into ${\rm int}(D_\varepsilon)$ (resp.
${\rm int}(-D_\varepsilon)$), then (iv) holds by the proof of
Proposition 3.2(ii) of \cite{LiuWaWe}. }
\end{remark}

We shall generalize Theorem~\ref{th:5.6} as follows.

\begin{theorem}\label{th:5.8}
Let $H$ be a Hilbert space, and let  $P_H\ne H$ be a closed cone,
i.e. $P_H=\bar P_H$ is convex and satisfies $\R^+\cdot P_H\subset
P_H$, $P_H\cap(-P_H)=\{\theta\}$. Suppose that  a $C^1$-functional
   ${\cal L}:H\to\R$ has critical point $\theta$ and satisfies the following properties.
\begin{description}
\item[(i)] ${\cal L}$ is bounded from below, and satisfies the (PS)
condition.

\item[(ii)] There exists a positive element $e\in P_H\setminus\{\theta\}$ such that the cone
$$
\hbox{$D:=\{u\in H\,|\, u\ge e\}\subset P_H$
 (resp. $-D$)}
$$
 contains all positive (resp. negative) critical
points of ${\cal L}$.

\item[(iii)] There exists a $\varepsilon_0\in (0, \|e\|)$ such that
for each $\varepsilon\in (0, \varepsilon_0]$,
\begin{equation}\label{e:5.9}
 \left.\begin{array}{ll}
 \hbox{${\cal K}:=id-\nabla{\cal L}: H\to H$ maps $D_\varepsilon$ (resp.
$-D_\varepsilon$)}\\
\qquad\hbox{into ${\rm int}(D_\varepsilon)$ (resp. ${\rm
int}(-D_\varepsilon)$)}.\end{array}\right\}
\end{equation}
\end{description}
Then there exists $\varepsilon_1\in (0, \varepsilon_0]$ such that
for any $0<\varepsilon\le\varepsilon_1$, $c_1$ defined by
(\ref{e:5.6}) is finite and $K^\ast_{c_1}:=\{u\in H\setminus
W_\varepsilon,|\, {\cal L}(u) = c_1,\;{\cal L}'(u)=0\}\ne\emptyset$.
Moreover, if $c_1<{\cal L}(\theta)$ and the conditions of
Theorem~2.1 of \cite{Lu1} are satisfied near each $u\in
K^\ast_{c_1}$ and $\theta$ then ${\cal L}$ has at least two
nontrivial critical points in $H\setminus W_\varepsilon$ provided
that (\ref{e:5.7}) and (\ref{e:5.8}) hold.
\end{theorem}

Since $D_\varepsilon$ and $-D_\varepsilon$ are positive invariant
under the pseudo-gradient flow $\widetilde\varphi^t$ defined by
(\ref{e:5.11}) (see proof below (\ref{e:5.12})), as in
Remark~\ref{rm:5.7}(i) we may show that ${\cal L}$ has a critical
point $u\in D_\varepsilon$ (resp. $u\in -D_\varepsilon$) with ${\cal
L}(u)=\inf_{v\in D_\varepsilon}{\cal L}(v)$ (resp. ${\cal
L}(u)=\inf_{v\in -D_\varepsilon}{\cal L}(v)$).
Remark~\ref{rm:5.7}(iii) shows that Theorem~\ref{th:5.8}(iii) is
stronger than Theorem~\ref{th:5.6}(iv). \vspace{2mm}

\noindent{\bf Proof of Theorem~\ref{th:5.8}.} The basic ideas is
almost the same as in \cite{LiuWaWe}. However, since we only assume
${\cal L}$ to be $C^1$, the negative gradient flow of it cannot be
used. We shall overcome this difficulty by some methods in
\cite{Ba2}.

Note that $V:=\nabla{\cal L}$ is a $C^0$ pseudo-gradient vector
field for ${\cal L}$ in the sense of \cite{Ba2}, i.e. $V:H\to H$ is
a continuous map satisfying
$$
\|V(x)\|<2\|\nabla{\cal L}(x)\|\quad\hbox{and}\quad (\nabla {\cal
L}(x), V(x))_H>\frac{1}{2}\|\nabla{\cal L}(x)\|^2\quad\forall x\in
H\setminus K({\cal L}).
$$
Moreover, (\ref{e:5.9}) means that $D$ and $-D$ are ${\cal
K}$-attractive in the sense of \cite[Definition 3.3]{Ba2}.  By
\cite[Lemma 3.4]{Ba2} there exists $\varepsilon_1\in (0,
\varepsilon_0]$ such that for every $\sigma\in (0, \varepsilon_1)$
there is a pseudo-gradient vector field $\widetilde V_\sigma$ of
${\cal L}$ such that for all $\varepsilon\in [\sigma,
\varepsilon_1]$ the sets $D_\varepsilon$ and $-D_\varepsilon$ are
strongly positive invariant under $-\widetilde V_\sigma$ in the
following sense:
\begin{equation}\label{e:5.10}
 \left.\begin{array}{ll}
 \hbox{for any $u\in\partial D_\varepsilon$, $\exists\,\epsilon_0>0$ such
 that $\forall\epsilon\in (0, \epsilon_0]$,}\\
 \hbox{$u+\epsilon \widetilde V_\sigma(u)\in{\rm Int}(D_\varepsilon)$ and $-u+\epsilon
\widetilde V_\sigma(-u)\in{\rm
Int}(-D_\varepsilon)$}.\end{array}\right\}
\end{equation}
Let $\widetilde {\cal K}_\sigma:=id-\widetilde V_\sigma$ and let
$\widetilde\varphi^t$ be the flow of $-\widetilde V_\sigma$, i.e.
\begin{equation}\label{e:5.11}
 \frac{d}{dt}\widetilde\varphi^t(u)=- \widetilde V_\sigma(\widetilde\varphi^t(u))\quad\hbox{and}\quad \widetilde\varphi^0(u)=u\in H.
\end{equation}
We claim:
\begin{equation}\label{e:5.12}
\left.\begin{array}{ll} \hbox{for any $\varepsilon\in
[\sigma,\varepsilon_0]$ the sets $D_\varepsilon$ and
$-D_\varepsilon$ are}\\
\hbox{strictly positive invariant for the flow
$\widetilde\varphi^t$.}\end{array}\right\}
\end{equation}

Indeed, from (\ref{e:5.10}) and Theorem 5.2 in \cite{Dei} we deduce
that the sets $D_\varepsilon$ and $-D_\varepsilon$ are positive
invariant for the flow $\widetilde\varphi^t$ for any $\varepsilon\in
[\sigma,\varepsilon_0]$.

As in the proof of \cite[Proposition 3.2(ii)]{LiuWaWe}, suppose by
contradiction that there exist $u\in D_\varepsilon$ and $t>0$ such
that $\widetilde\varphi^t(u)\in\partial D_\varepsilon$. Then Mazur's
separation theorem yields $f\in H^\ast$ and $\beta>0$ such that
$f(\widetilde\varphi^t(u))=\beta$ and $f(u)>\beta$ for any $u\in{\rm
Int}(D_\varepsilon)$. By (\ref{e:5.10}) we have $\epsilon_0>0$ such
that
 $\widetilde\varphi^t(u) +\epsilon \widetilde
V(\widetilde\varphi^t(u))\in{\rm
Int}(D_\varepsilon)\;\forall\epsilon\in (0, \epsilon_0]$. Hence
\begin{eqnarray*}
\frac{d}{ds}f(\widetilde\varphi^s(u))\Bigm|_{s=t}&=&f\bigl(-\widetilde
V(\widetilde\varphi^t(u))\bigr)=\frac{1}{\epsilon}f\bigl(-\epsilon\widetilde
V(\widetilde\varphi^t(u))\bigr)\\
&=&\frac{1}{\epsilon}f\bigl(\widetilde\varphi^t(u)-\epsilon\widetilde
V(\widetilde\varphi^t(u))\bigr)-\frac{1}{\epsilon}f\bigl(\widetilde\varphi^t(u)\bigr)\\
&=&\frac{1}{\epsilon}f\bigl(\widetilde\varphi^t(u)-\epsilon\widetilde
V(\widetilde\varphi^t(u))\bigr)-\frac{1}{\epsilon}\beta>0.
\end{eqnarray*}
This leads to a contradiction as in \cite{LiuWaWe}. The assertion in
(\ref{e:5.11}) is proved.

Having (\ref{e:5.12}) we may follow the lines of \cite{LiuWaWe} to
outline the remaining proof.

 Let $\bar H^\ast$ denote Alexander-Spanier cohomology with
coefficients in the field $\Z_2$.  For a critical point $u$ of
${\cal L}$ the cohomological critical groups of ${\cal L}$ at $u$
are defined by
$$
C^q({\cal L}, u):=\bar H^q({\cal L}^c, {\cal
L}^c\setminus\{u\})\quad\forall q\ge 0,
$$
 where $c={\cal L}(u)$. We conclude
\begin{equation}\label{e:5.13}
C^q({\cal L}, u)\cong C_q({\cal L}, u;\Z_2)\quad\forall q\in\Z.
\end{equation}
This can be obtained from the proof of Proposition~\ref{prop:3.3}.
We can also prove it as follows.  By excision property it holds that
$C^q({\cal L}, u)=\bar H^q({\cal L}^c\cap U, ({\cal
L}^c\setminus\{u\})\cap U)\;\forall q\ge 0$ for any open
neighborhood $U$ of $u$. In particular, one can find a $U$ so that
both ${\cal L}^c\cap U$ and $({\cal L}^c\setminus\{u\})\cap U$ are
absolute neighborhood retracts ( \cite[Th.1.1]{Deg} and \cite[Remark
2]{Cor1}). It follows that
$$
\bar H^q({\cal L}^c\cap U, ({\cal L}^c\setminus\{u\})\cap U)\cong
H^q({\cal L}^c\cap U, ({\cal L}^c\setminus\{u\})\cap U)\;\forall
q\ge 0.
$$
Moreover, $H^q({\cal L}^c\cap U, ({\cal L}^c\setminus\{u\})\cap U;
G)\cong{\rm Hom}(H_q({\cal L}^c\cap U, ({\cal
L}^c\setminus\{u\})\cap U; G), G)$ for any divisible group $G$.
 If $C^1({\cal L}, u)\ne 0$ we say $u$ to be of {\it
mountain-pass type}.

\begin{lemma}\label{lem:5.9}{\rm (\cite[Lemma 4.4]{LiuWaWe})}. If $C^1({\cal L}, u)\ne
0$, and $\nu(u)\le 1$ in case $\mu(u)=0$, then $C^q({\cal L},
u)=\delta_{q1}\Z_2$ for $q\in\Z$.
\end{lemma}

\noindent{\bf Proof}. Note that $\bar H^q(B^n, S^{n-1};
G)=\delta_{qn}G$. If $\nu(u)={\rm Ker}({\cal L}''(u))=0$, i.e. $u$
is nondegenerate, by Morse lemma we get that $C^q({\cal L}, u)\cong
H^q(B^{\mu(u)}, S^{\mu(u)-1};\Z_2)=\delta_{q\mu(u)}\Z_2$. The
desired conclusion follows.

 If $\nu(u)={\rm Ker}({\cal L}''(u))>0$, from Corollary~2.6 of \cite{Lu1} it follows that
$\mu(u)\le 1\le\mu(u)+\nu(u)$ and $C^1({\cal L}, u)\cong
C^{1-\mu(u)}({\cal L}^\circ, \theta)$, where ${\cal L}^\circ$ is a
function defined near the origin of a $\nu(u)$-dimensional space
${\rm Ker}({\cal L}''(u))$. If $\mu(u)=1$ then Corollary~2.9(iii) of
\cite{Lu1} gives the conclusion. If $\mu(u)=0$ then $\nu(u)=1$ by
the fact that $1\le\mu(u)+\nu(u)=\nu(u)$. Corollary~2.9(ii) of
\cite{Lu1} yields the conclusion. $\Box$\vspace{2mm}

The following is a special case of the strong excision property
(Theorem 5) on the page 318 of \cite{Spa}.

\begin{lemma}\label{lem:5.10}{\rm (\cite[Lemma 4.4]{LiuWaWe})}. Let $X$
be paracompact Hausdorff space and let $Y, Z\subset
X$ be closed subsets such that $X=Y\cup Z$. Then the inclusion $(Y,
Y\cap Z)\to (X, Z)$ induces an isomorphism $\bar H^\ast(X, Z)\to\bar
H^\ast(Y, Y\cap Z)$.
\end{lemma}

By the assumption (i) in Theorem~\ref{e:5.8}, ${\cal L}$ is bounded
from below, for each $u\in H$ the flow
$t\mapsto\widetilde\varphi^t(u)$ exists in an open interval
containing $[0, \infty)$. For $-\infty\le c\le \infty$, ${\cal L}^c$
is positively invariant for the flow $\widetilde\varphi^t$, and
strictly positively invariant if $c$ is a regular value of ${\cal
L}$.

\begin{lemma}\label{lem:5.11}{\rm (\cite[Lemma 4.3]{LiuWaWe})}.
  Let $u$ be a critical point of ${\cal L}$ with ${\cal L}(u)=c$
  and such that $B_H(u, 2\varepsilon)$ contains no other critical point of ${\cal L}$.
Then for $\delta>0$ sufficiently small there exists a closed
neighborhood $N\subset B_H(u,\varepsilon)\cap {\cal L}^{c+\delta}$
of of $u$ such that $N\cup{\cal L}^{c-\delta}$  is positively
invariant for the flow $\widetilde\varphi^t$. Moreover,
$$
C^q({\cal L},u)\cong \bar H^q(N\cup{\cal L}^{c-\delta}, {\cal
L}^{c-\delta})\cong \bar H^q(N, {\cal L}^{c-\delta}\cap
N)\quad\forall q\in\Z.
$$
\end{lemma}

Fix $\varepsilon\in [\sigma, \varepsilon_1]$, then
$W_\varepsilon:=D_\varepsilon\cup(-D_{\varepsilon})$ is closed and
strictly positively invariant under $\widetilde\varphi^t$ by
(\ref{e:5.11}). In particular, $\partial W_\varepsilon$ contains no
critical points of ${\cal L}$. So the (PS) condition implies that
$K^\ast_c:=\{u\in H\setminus W_\varepsilon,|\, {\cal L}(u) =
c,\;{\cal L}'(u)=0\}$
  is a compact subset of $H\setminus W_\varepsilon$ for every $c\in\R$.
Clearly, every non-trivial critical point in $K^\ast_c$ changes
sign.

\begin{lemma}\label{lem:5.12}{\rm (\cite[Lemma 4.5]{LiuWaWe})}.
Suppose that $-\infty\le a<b\le c\le d\le\infty$ satisfy:
\begin{description}
\item[(i)] $K^\ast_{c'}=\emptyset$ for any $c'\in [b, d]\setminus\{c\}$

\item[(ii)] There is a neighborhood $N\subset{\cal L}^d$ of $K^\ast_c$ such that ${\cal L}^b\cup N$ is
is positively invariant under $\widetilde\varphi^t$.
\end{description}
 Then the inclusion
 $$
 (N\cup{\cal L}^b\cup W_\varepsilon, {\cal L}^a\cup W_\varepsilon)
\to ({\cal L}^d\cup W_\varepsilon, {\cal L}^a\cup W_\varepsilon)
$$
is a homotopy equivalence. In particular, if $K^\ast_c=\emptyset$
for any $c\in [b, d]$, then the inclusion
 $({\cal L}^b\cup W_\varepsilon, {\cal L}^a\cup W_\varepsilon)
\to ({\cal L}^d\cup W_\varepsilon, {\cal L}^a\cup W_\varepsilon)$ is
a homotopy equivalence.
\end{lemma}

From Lemmas~\ref{lem:5.10}, \ref{lem:5.11} and \ref{lem:5.12} we
have

\begin{lemma}\label{lem:5.13}{\rm (\cite[Lemma 4.7]{LiuWaWe})}.
If the compact set $K^\ast_c$ consists of isolated critical points
$u_1, \cdots, u_m$ for some $c\in\R$, then
 $$
 \bar H^q({\cal L}^{c+\delta}\cup W_\varepsilon, {\cal L}^{c-\delta}\cup
 W_\varepsilon)\cong \bigoplus^m_{i=1} C^q({\cal L}, u_i)
 $$
for $q\in\Z$ and $\delta>0$ sufficiently small.
\end{lemma}

Note that
\begin{equation}\label{e:5.14}
\bar H^{q}(H, W_\varepsilon;\Z_2)\cong {\rm Hom}(H_{q}(H,
W_\varepsilon;\Z_2);\Z_2)\cong\delta_{q1}\Z_2.
\end{equation}
By the definition of $c_1$ and the assumption (i) we get
$c_1>-\infty$.

\begin{lemma}\label{lem:5.14}{\rm (\cite[Lemma 4.8]{LiuWaWe})}.
\begin{description}
\item[(i)] $K^\ast_{c_1}\ne\emptyset$, and $c_1\le {\cal L}(\theta)$ provided that
 ${\cal L}(te)\le{\cal L}(\theta)\;\forall t\in [-1,
1]$.

\item[(ii)] If $c_1<{\cal L}(\theta)$ then $K^\ast_{c_1}$ consists of
non-trivial critical points.

\item[(iii)] If  $K^\ast_{c_1}$ consists of only one isolated critical point
$u$, then
 $$
C^q({\cal L}, u)\cong \bar H^q({\cal L}^{c_1+\delta}\cup
W_\varepsilon, {\cal L}^{c_1-\delta}\cup
 W_\varepsilon)\cong \delta_{q1}\Z_2\quad\hbox{for small
 $\delta>0$}.
 $$
 \end{description}
\end{lemma}

The claim that $c_1\le{\cal L}(\theta)$ in (i) can be proved by
(\ref{e:5.14}). Lemma~\ref{lem:5.12} leads to
$K^\ast_{c_1}\ne\emptyset$. Then the assumption and
Lemma~\ref{lem:5.13} imply $C^q({\cal L}, u)\cong\bar H^q({\cal
L}^{c+\delta}\cup W_\varepsilon, {\cal L}^{c-\delta}\cup
 W_\varepsilon)$ for small $\delta>0$. From this and the definition
 of $c_1$ it follows that $C^1({\cal L}, u)\ne 0$. Combing it with
 (\ref{e:5.7}) together, we may use Lemma~\ref{lem:5.9} to infer (iii).

Now (\ref{e:5.8}), (\ref{e:5.14}) and
Lemmas~\ref{lem:5.11},~\ref{lem:5.12} and ~\ref{lem:5.14} lead to

\begin{lemma}\label{lem:5.15}{\rm (\cite[proposition 4.9]{LiuWaWe})}.
Under the assumption (\ref{e:5.8}), if $c_1<{\cal L}(\theta)$ and
then there exists  $c\ne c_1$ such that $K^\ast_{c}$ contains a
non-trivial  critical point of ${\cal L}$.
\end{lemma}

Summarizing these two lemmas  we complete the proof for
$\varepsilon\in [\sigma,\varepsilon_1]$. Since $\sigma\in
(0,\varepsilon_0)$ is arbitrary Theorem~\ref{th:5.8} is proved.
$\Box$\vspace{2mm}

\part{Applications }\label{p:2}

In progress!

\end{document}